\documentclass[11pt,a4paper,leqno]{article}
\usepackage{amsmath}
\usepackage{amssymb}
\usepackage{graphicx}
\advance\topmargin by -2.2cm
\newcommand{\bbr}{I\!\!R}

\newcommand{\XXX}{{\bar X}}

\newcommand{\barr}{\begin{array}}
\newcommand{\earr}{\end{array}}
\newcommand{\beqq}{\begin{equation}}
\newcommand{\eeqq}{\end{equation}}
\newcommand{\beao}{\begin{eqnarray*}}
\newcommand{\eeao}{\end{eqnarray*}\noindent}
\newcommand{\beam}{\begin{eqnarray}}
\newcommand{\eeam}{\end{eqnarray}\noindent}

\newcommand{\halmos}{\quad\hfill\mbox{$\Box$}}

\newcommand{\al}{\alpha}

\newcommand{\vep}{\varepsilon}

\newtheorem{theo}{Theorem}
\newtheorem{prop}{\indent Proposition}

\newtheorem{defin}{\indent Definition}
\newtheorem{cor}{\indent Corollary}

\newcommand{\wt}{\widetilde}

\newcommand{\ov}{\overline} 
\newcommand{\ul}{\underline}

\author{R. H\"opfner$^*$ \and E.~L\"ocherbach \and M. Thieullen\thanks{This work has been supported by the Agence Nationale de la Recherche through the project MANDy, Mathematical Analysis of Neuronal Dynamics, ANR-09-BLAN-0008-01. e-mail addresses: \tt{hoepfner@mathematik.uni-mainz.de},
     \tt{eva.loecherbach@u-cergy.fr} and
    \tt{michele.thieullen@upmc.fr}} \\
{\it Johannes Gutenberg-Universit\"at Mainz, Universit\'e de Cergy-Pontoise} \\
 {\it and Universit\'e Pierre et Marie Curie.}}
\begin{document}

\title{Transition densities for strongly degenerate time inhomogeneous random models.}

\maketitle

\begin{abstract}
In this paper we study the existence of densities for strongly degenerate stochastic differential equations whose coefficients depend on time and are not globally Lipschitz. In these models neither local ellipticity nor the strong H\"ormander condition is satisfied. In this general setting we show that continuous transition densities indeed exist in all neighborhoods of points where the weak H\"ormander condition is satisfied. We also exhibit regions where these densities remain positive. We then apply these results to stochastic Hodgkin-Huxley models as a first step towards the study of ergodicity properties of such systems.

\end{abstract}

\\


\section{Introduction}
In this paper we study the existence of densities for strongly degenerate stochastic differential equations (SDEs) whose coefficients depend on time and are not globally Lipschitz. Such models are not covered by previous results in the literature. We consider multidimensional stochastic systems where the noise is one dimensional and present only in the first and last component. Neither local ellipticity nor the strong H\"ormander condition is satisfied. 

Our original motivation was to be able to describe with probabilistic tools the long time behavior of a neuron embedded in a network, receiving synaptic stimulation from a large number of other neurons in the network  through its dendritic tree.  More precisely this stimulation takes the form of a random input carrying a deterministic and periodic signal, and the question is the ergodicity of the process resulting in the neuron on the one hand and the input it receives on the other hand. In this model, the neuron part is based on the Hodgkin-Huxley model well-known in physiology and the signal part is a noisy perturbation of a deterministic signal. This particular question leads us to the study of a non time homogeneous 5-dimensional stochastic system driven by a one dimensional Brownian noise. Prior to adressing questions such as ergodicity properties for such systems, e.g. in the sense of \cite{MT1}-\cite{MT2}, the first step is to establish that Lebesgue densities exist and --at least on suitable parts of the state space-- are continuous and strictly positive. This is the topic of the present paper where we study densities for a general class containing in particular the stochastic conductance-based models well known in physiology and among these the stochastic Hodgkin-Huxley system. In a companion paper \cite{HLT-2} we address the periodic ergodicity of the Hodgkin-Huxley model when driven by an Ornstein-Uhlenbeck perturbation of a deterministic signal.

In order to study the existence of densities and their regularity it is now classical to make use of Malliavin calculus techniques. In this theory the H\"ormander condition plays a fundamental role as a sufficient condition for the existence of densities. This condition is satisfied if the Lie algebra generated by the coefficients of the SDE has full dimension which means in other words that the diffusion is actually strong enough even if the noise is visible only on a restricted number of components. The H\"ormander condition has two forms: the strong form involves only Lie brackets computed using the diffusion coefficients, whereas the weak form may also include the drift coefficient. For non time homogeneous SDEs, to the best of our knowledge, the existing results all require at least the strong H\"ormander condition (see \cite{C-M} and the references therein). 

For our models, we can only hope to satisfy the weak H\"ormander condition since the diffusion vector field is one dimensional. In general this condition will hold only locally. This is well illustrated in the last section when we address the stochastic Hodgkin-Huxley model. SDEs satisfying local H\"ormander condition with locally smooth coefficients have been considered recently in a time homogeneous setting (cf. \cite{Bally}, \cite{stefano}, \cite{Fournier}) where the local H\"ormander condition is ensured by a local ellipticity assumption so it is a strong H\"ormander condition that is considered in these works.  We extend these results to our framework using a technique based on estimates of the Fourier transform introduced in these papers. More precisely, we use a localization argument which is based on ideas of \cite{stefano}. However our frame is technically more difficult since time homogeneity fails and the H\"ormander condition, which holds only locally, is the weak one. In this general setting we show that continuous transition densities indeed exist in all neighborhoods of points where the weak H\"ormander condition is satisfied. We also prove that these densities are lower semi continuous (lsc) w.r.t. the starting point even if our system does not enjoy the Feller property.

We are then able to say more when we consider a family of SDEs that we call {\it SDEs with internal variables and random input.} They are of the following form
\begin{eqnarray}\label{internalintro}
dX_{1,t}&=& F (X_{l,t}, \ 1\leq l\leq m-1)dt +dX_{m,t} \, ,\\
dX_{i,t}&=&[-a_i(X_{1,t})X_{i,t}+b_i(X_{1,t})]dt, \, \, \, i = 2, \cdot, \cdot, \cdot , m-1\, ,\nonumber\\
dX_{m,t} &=& b_m(t,X_{m,t}) dt \;+\sigma (X_{m,t})dW_t.\nonumber 
 \end{eqnarray}
In (\ref{internalintro}) the last component $X_m$ can be seen as an external random input acting on the deterministic evolution
 \begin{eqnarray}\label{subinternalintro}
dz_{1,t}&=& F (z_{l,t}, \ 1\leq l\leq m-1) dt  \, ,\\
dz_{i,t}&=&[-a_i(z_{1,t})z_{i,t}+b_i(z_{1,t})]dt, \, \, \, i = 2, \cdot, \cdot, \cdot , m-1.\nonumber
 \end{eqnarray}

These systems have an important interpretation related to modeling in physics and biology. Indeed one can show that (\ref{subinternalintro}) is the limit of a sequence of stochastic processes in the sense of the Law of Large Numbers or Fluid Limit (cf. \cite{CDMR}, \cite{RGGMK}, \cite{FGRC}, \cite{PTW-10}). Each stochastic process in the sequence is a Piecewise Deterministic Markov Process (cf. \cite{D}). More precisely the interpretation of (\ref{subinternalintro}) is the following. Consider a population of individuals of different types $ i\in \{2,\cdot, \cdot, m-1\}$, each individual of type $i$ being in two states, active or inactive (open or closed). The individuals are coupled by a global variable $z_1$ since their transition rates $a_i(z_1)$, $b_i(z_1)$ depend on $z_1$. If there are $N$ individuals of each type, consider at any time $t$, the vector consisting in $z_1$ and the proportions of individuals of each type in the active state. This stochastic vector is a Piecewise Deterministic Markov Process. When $N\mapsto \infty$, the limit of the sequence is (\ref{subinternalintro}). 
Then $z_{i,t}$ gives the probability that an individual of type $i$ is active at time $t$. The detailed form of these fonctions in the application we have in mind will be provided in Section \ref{sec:main}. In this particular model derived from the model of Hodgkin and Huxley (cf. \cite{HH-52}), $z_1$ describes the evolution of the membrane potential of a neuron, which can be observed. This evolution results from the gating mechanism of the ion channels located in the membrane $z_i$, $i\in\{2,3,4\}$ that are not observed. In that model, $m=5$. Note that the whole class of conductance-based models classical in physiology can be included in this setting. The Hodgkin-Huxley model belongs to this class. Another example is the Morris-Lecar model for the excitation of muscle fibers first presented in \cite{ML-81}. For conductance-based models it is important to note that $F (x_l, \ 1\leq l\leq m-1)$ in (\ref{internalintro}) is linear w.r.t. $x_1$ as a consequence of Ohm's Law. We give more details in section \ref{sec:main}.

For such systems we provide an explicit form of the weak H\"ormander condition via the computation of a determinant built up with the coefficients of (\ref{subinternalintro}) and their successive derivatives. Then we exhibit regions where densities, if they exist, remain positive. These regions are related to neighborhoods of equilibrium points of (\ref{subinternalintro}). To prove this result we exploit the particular structure of (\ref{subinternalintro}) namely the linearity of the internal equations (those for $dz_i, \, i = 2, \cdot, \cdot, \cdot , m-1$) which ensures global asymptotic stability of the corresponding equilibria when $z_1$ is held fixed. 

We also show that with positive probability, the solution of (\ref{internalintro}) can imitate any deterministic evolution resulting from an arbitrary input applied to (\ref{subinternalintro}), on an arbitrary interval of time.

In the last section we apply the previous results to stochastic conductance-based models and in particular Hodgkin-Huxley systems. As mentioned previously they belong to the family of SDEs with internal variables and random input. For the Hodgkin-Huxley model we conduct a numerical study of the determinant which shows that the H\"ormander condition is satisfied at certain stable equilibrium points and/or along a specific stable periodic orbit. Therefore, depending on the starting point, the stochastic Hodgkin-Huxley model possesses strictly positive densities either in small neighborhoods of the above mentioned equilibrium point or of the periodic orbit.Ê

The present paper is the first of two papers in which we propose a probabilistic study of the periodic ergodicity for stochastic Hodgkin-Huxley models. To the best of our knowledge, no other probabilistic study has been presented in the literature before. There are some simulation studies (see. e.g. \cite{PPM-05} and \cite{YWWL-01}), but not much seems to be known mathematically. In this first part of our study we look for densities of the system. We refer the reader to the companion paper \cite{HLT-2} for the following step showing recurrence properties when the stochastic input is a mean reverting process of Ornstein-Uhlenbeck type carrying a deterministic periodic signal.Ê

Our paper is organized as follows. In section \ref{models} we first present our general models and assumptions as well as the family of SDEs with internal variables and random input. Section \ref{density general} is devoted to the proof of the existence of densities locally for time inhomogeneous systems with locally Lipschitz coefficients. In section \ref{ivri} we explicit the weak H\"ormander condition and address the positivity of densities for SDEs with internal variables and random input. The last section of the paper is devoted to the deterministic and stochastic Hodgkin-Huxley models. We apply the results obtained in the previous sections, in particular through a numerical study. We end up this section with an application to conductance-based models.

\section{Our Models}\label{models}

\subsection{General Assumptions} 

We now describe the general framework of our results. Given an integer $ m \geq 1 ,$ we consider processes taking values in $\bbr^m$ and write $ x = (x_1, \ldots , x_m) $ for generic elements of $\bbr^m .$ Let $\sigma $ be a measurable function from $\bbr^m $ to $\bbr^m $ and $b$ a smooth function from $ [0, \infty [ \times \bbr^m $ to $\bbr^m . $ For all $ x \in \bbr^m,$ we consider the SDE 
\begin{equation}\label{eq:equationms}
X_{i,t}= x_i + \int_0^t b_i ( s, X_s) ds + \int_0^t \sigma_i (X_s) d W_s , \; t \geq 0 , \; i = 1, \ldots , m. 
\end{equation}
Here, $W$ is a one-dimensional Brownian motion and $ \sigma$ is identified with an $m \otimes 1-$matrix. We write $P_x$ for the probability measure under which the solution $X = (X_t)_{t \geq 0 } $ of (\ref{eq:equationms}) starts at $x.$ 

We work under the assumptions {\bf (H1)}, {\bf (H2)}:

{\bf (H1)} Equation (\ref{eq:equationms}) admits unique strong solutions. Whenever $ ( X_t, t \geq 0 )$ is such a strong solution with $X_0=x$, we assume that 
there exists an increasing sequence of compacts $ \; K_n \subset K_{n+1}$ such that, if $x$ belongs to
$\bigcup_{n } K_n ,$ then $T_n := \inf \{ t : X_t \notin K_n \}  \to \infty  \mbox{ almost surely as } n \to \infty .$ We take $K_n = [a_n, b_n]= \prod_{i = 1}^m [a_{n,i} , b_{n,i} ] ,$ where $ a_n = ( a_{n,1} , \ldots , a_{n,m} )$. 
Due to the latter assumption, the state space of the process $ ( X_t, t \geq 0 )$ is $E := \bigcup_{n } K_n.$

{\bf (H2)} The coefficients of  (\ref{eq:equationms}) are locally smooth. Namely we assume that for all $n$, $\, \sigma \in C_b^\infty (K_n , \bbr^m )$ and for every multi-index $\beta \in \{ 0, \ldots, m \}^l, l \geq 1, \quad  b(t,x) +  \partial_{\beta}  b (t,x)  \mbox{ is bounded on  }  [ 0, T ] \times K_n$ for all $T>0$. Here $\partial_\beta =\frac{ \partial^l}{ \partial x_{\beta_1} \ldots \partial x_{\beta_l} } $ and we identify $x_0 $ with $t.$

\subsection{Towards Application: internal variables and random input}\label{sec:internal}

As mentioned in the introduction, our original motivation in the present paper is to investigate the long time behavior of a neuron in a network receiving a dendritic input, using a probabilistic approach. This leads us to consider a subclass of models where the SDE (\ref{eq:equationms}) possesses a particular structure as follows

\begin{eqnarray}\label{internal1}
dX_{1,t}&=& F (X_{l,t}, \ 1\leq l\leq m-1)dt +dX_{m,t}\, , \\
dX_{i,t}&=&[-a_i(X_{1,t})X_{i,t}+b_i(X_{1,t})]dt, \, \, \, i = 2, \cdot, \cdot, \cdot , m-1\, ,\nonumber\\
dX_{m,t} &=& b_m(t,X_{m,t}) dt \;+\sigma (X_{m,t})dW_t\, ,\nonumber 
 \end{eqnarray}
for some functions $F$, $a_i$, $b_i$, $b_m$, $\sigma$. We will work under the following assumption:

\noindent {\bf (H3)} The SDE 
\begin{equation*}
dZ_t = b_m(t,Z_t) dt \;+\sigma (Z_t)dW_t
\end{equation*}
possesses a unique strong solution $(Z_t)_{t\ge 0}$ taking values in an open interval $U\subset \bbr$. Moreover $\sigma(\cdot)$ is strictly positive on $U$ and its restriction to every compact interval in $U$ is of class $C^\infty$.  

The random signal $X_m$ is an external input to the underlying deterministic system
\begin{eqnarray}\label{subinternal}
dz_{1,t}&=& F (z_{l,t}, \ 1\leq l\leq m-1) dt \, ,  \\
dz_{i,t}&=&[-a_i(z_{1,t})z_{i,t}+b_i(z_{1,t})]dt, \, \, \, i = 2, \cdot, \cdot, \cdot , m-1.\nonumber
 \end{eqnarray}
This system can be interpreted as the limit of a sequence of stochastic ones in the sense of the Law of Large Numbers or Fluid Limit (cf \cite{PTW-10}): consider a population of individuals of different types $ i\in \{2,\cdot, \cdot, m-1\}$, each individual of type $i$ being in two states, active or inactive (open or closed). The individuals are coupled by a global variable $z_1$, their transition rates $a_i(z_1)$, $b_i(z_1)$ depending on $z_1$. If there are $N$ individuals of each type, consider at any time $t$, the vector consisting in $z_1$ and the proportions of individuals of each type in the active state. When $N\mapsto \infty$, the limit of this process is (\ref{subinternal}). Then $z_{i,t}$ gives the probability that an individual of type $i$ is active at time $t$. The detailed form of these fonctions in the application we have in mind will be provided in Section \ref{sec:main}. In this particular model derived from the model of Hodgkin and Huxley (cf. \cite{HH-52}), $z_1$ describes the evolution of the membrane potential of a neuron, which can be observed. This evolution results from the gating mechanism of specific ion channels located in the membrane $z_i$, $i\in\{2,3,4\}$ that are not observed. 
 
The linear form of the equations for the variables $X_i$ , $i\notin \{1,m\}\, $ in (\ref{internal1})  has important consequences as we recall in the following proposition. 
\begin{prop}\label{linear}
Fix $i\in\{2,\cdot\cdot\cdot,m-1\}$. Assume that $X_{i,0}\in[0,1]$ a.s., and on the interval $[0,1]$ the function $a_i$ is positive and $b_i(\cdot)\leq a_i(\cdot)$. Then $\forall t>0,\, X_{i,t}\in[0,1]$ a.s.
\end{prop} 
 
{\bf Proof.} Given the trajectory of $X_1 $, the variation of constants method yields the following representation
\begin{equation}\label{eq:exact}
X_{i,t}=X_{i,0}{\rm e}^{-\int_0^t a_i(X_{1,s})ds}+\int_0^t b_i(X_{1,u}){\rm e}^{-\int_u^t a_i(X_{1,r})dr}du.
\end{equation}
However note that (\ref{eq:exact}) does not provide an explicit formula for $X_{i,t}$ since $X_1 $ depends on $X_i$ (the system is fully coupled).  

Writing $\int_0^t b_i(X_{1,u}){\rm e}^{-\int_u^t a_i(X_{1,r})dr}du=\int_0^t \frac{b_i(X_{1,u})}{a_i(X_{1,u})}a_i(X_{1,u}){\rm e}^{-\int_u^t a_i(X_{1,u})dr}du$, the assumptions on $a_i(\cdot)$ and $b_i(\cdot)$ imply 
\begin{equation}
0 \leq   X_{i,t}\leq X_{i,0}{\rm e}^{-\int_0^t a_i(X_{1,s})ds}+\int_0^t a_i(X_{1,u}){\rm e}^{-\int_u^t a_i(X_{1,r})dr}du.
\end{equation}
By straightforward integration it follows that
\begin{eqnarray}\label{eq:tobecited}
0 \leq  X_{i,t}&\leq& (X_{i,0}+{\rm e}^{\int_0^t a_i(X_{1,r})dr}-1){\rm e}^{-\int_0^t a_i(X_{1,s})ds} \nonumber \\
&=& 1+(X_{i,0}-1){\rm e}^{-\int_0^t a_i(X_{1,s})ds}.
\end{eqnarray}
The statement follows. \halmos

\begin{prop}\label{globallyasymptoticallystable}
Take two real numbers $a$ and $b$, with $a>0$ and define $y_\infty:=\frac{b}{a}$. Denote by $y_t$ the solution of the ode $dy_t = (-a y_t+b)dt$. For all $\varepsilon>0$ there exists $ t_0$ such that 
\begin{equation*}
|y_t-y_\infty|<\varepsilon\quad\quad\forall t\geq t_0.
\end{equation*}
\end{prop}

\noindent{\bf Proof} This well known result follows from the positivity of $a$ and the explicit form of $y_t$ given by $y_t =( y_0 -y_\infty)e^{ - at} +y_\infty$. \halmos

\noindent In particular, if in (\ref{subinternal}) the $a_i$ are all positive and if $z_1$ were kept constant, the vector $(z_{i,t}, \,i = 2, \cdot, \cdot, \cdot , m-1)$ would converge when $t\rightarrow +\infty$ to its equilibrium $(y_{i,\infty}=\frac{b_i(z_1)}{a_i(z_1)},\, \,i = 2, \cdot, \cdot, \cdot , m-1)$ which is globally asymptotically stable. We will use the above result in Section \ref{positivedensities} below.

\section{Smoothness of densities of a strongly degenerate SDE with locally smooth coefficients depending on time}\label{density general}
Classically, one
proves that the solution of an SDE admits a smooth density via Malliavin calculus and H\"ormander condition. Many authors assume that the coefficients of the SDE are $C^\infty$, bounded, with bounded derivatives of any order and that H\"ormander condition is satisfied all over the state space. In the application we are interested in, the regularity condition on the coefficients is not satisfied since they are not globally Lipschitz. Regarding the H\"ormander conditions there are actually two possibilities: either to work under the strong H\"ormander condition or under the weak one which is a less stringent assumption. Many authors work under the strong H\"ormander condition whereas we work here under the weak one which moreover holds only locally. This is because the systems we consider are highly degenerate. In addition the drift coefficient depends on time. For these reasons, we have to apply local arguments in a non time homogeneous setting. In what follows we extend the results of \cite{KS} which hold only in a time homogeneous framework. Then, to prove local existence of densities, we rely on \cite{stefano} that we adapt to our framework.

\subsection{Local H\"ormander condition in a time dependent setting}\label{sec:hoerm} 
In this section we state our local weak H\"ormander condition. We first extend the results of \cite{KS} to time dependent coefficients.  The study of the H\"ormander condition requires to rewrite the SDE (\ref{eq:equationms})
in Stratonovitch form and see the obtained coefficients as linear differential operators of degree one (or vector fields). This amounts to replace the drift $ b(t,x) $ by 
$$ \tilde b_i (t,x) := b_i (t,x) - \frac12 \sum_{k=1}^m \sigma_k (x) \frac{ \partial \sigma_i }{\partial x_k } (x) , \; 1 \le i \le m \; , x \in \bbr^m\, , $$
which is again non-homogeneous in time and to introduce the vector field   
$$ A_0 = \frac{\partial }{\partial t } + \sum_{ i = 1}^m \tilde b_i (t,x) \frac{\partial }{\partial x_i } = \frac{\partial}{\partial t} + \tilde b $$
on $[0, \infty [ \times \bbr^m \to \bbr^m $. It can be identified with the $ (m+1 )-$dimensional function $ A_0 ( t, x ) = ( 1, \tilde b_1 , \ldots , \tilde b_m) .$ In what follows all functions $ {\cal T} (t,x) : [0, \infty [ \times \bbr^m \to \bbr^m $ different from $A_0$ will also be interpreted as vector fields
$$ {\cal T}(t,x) = \sum_{i=1}^m {\cal T}_i (t,x) \frac{\partial}{\partial x_i}$$
and identified with the $(m+1)-$dimensional function ${\cal T}(t,x) = (0, {\cal T}_1 , \ldots , {\cal T}_m ).$ 
The formalism in \cite{KS} uses multi-indices. Let $ M := \{ \emptyset \}\cup \bigcup_{ l = 1}^\infty  \{ 0 , 1 \}^l .$ For any $ \alpha \in M, $ define $ | \alpha | := l $ if $ \alpha \in \{ 0, 1\}^l, l \geq 1$ and $ | \emptyset | := 0 .$ Moreover, let $ \| \alpha \| := | \alpha | + card \{ j : \alpha_j = 0 \} $ if $ | \alpha | \geq 1 $, and $ \| \emptyset \| := 0. $ Finally, $ \alpha ' := ( \alpha_1, \ldots , \alpha_{l-1}) $ if $\alpha = ( \alpha_1 , \ldots , \alpha_l ), l \geq 2 , $ $\alpha' := \emptyset $ if $l = 1 .$ 

To build the successive Lie brackets we start with $ A_1 (x) := \sigma (x) $ that we identify with $ \sum_{ i =1}^m A_{1,i} (x) \frac{\partial }{ \partial x_i } \equiv \sum_{ i =1}^m \sigma_i (x) \frac{\partial }{ \partial x_i }$ and for $ {\cal T}  : [ 0, \infty [ \times \bbr^m \to \bbr^m ,$ we define $ {\cal T}_\emptyset (t,x) := {\cal T}(t,x) $ and for $| \alpha | \geq 1,$ $ {\cal T}_{(\alpha ) } ( t,x) := [ A_{ \alpha_l } , {\cal T}_{( \alpha ')} ] $ inductively in $|\alpha|.$ Here, $[{\cal T}, {\cal V}]$ denotes the Lie bracket defined by $ [{\cal T}, {\cal V} ]_i = \sum_{j = 0}^m \left( {\cal T}_j \frac{ \partial {\cal V}_i }{\partial x_j } - {\cal V}_j \frac{ \partial {\cal T}_i }{\partial x_j } \right) .$

In particular, if $ {\cal T} = A_1 \equiv \sigma $ we have $ [A_1, {\cal V}]_i = \sum_{j = 1}^m \left( \sigma_j \frac{ \partial {\cal V}_i }{\partial x_j } - {\cal V}_j \frac{ \partial \sigma_i }{\partial x_j } \right)$ and the time variable does not play any role. But if $ {\cal T} = A_0, $
$$ [A_0, {\cal V} ]_i = \sum_{j = 0}^m \left( A_{0,j} \frac{ \partial {\cal V}_i }{\partial x_j } - {\cal V}_j \frac{ \partial A_{0,i} }{\partial x_j } \right) = 
\frac{ \partial {\cal V}_i }{\partial t } +\sum_{j = 1}^m \left( A_{0,j} \frac{ \partial {\cal V}_i }{\partial x_j } - {\cal V}_j \frac{ \partial A_{0,i} }{\partial x_j } \right) .$$

We need additional notation in order to state the H\"ormander condition. For $y \in \bbr^m $ and $\delta > 0 $, we denote by $B_\delta (y) $ the open ball of radius $\delta $ centered at $y$. For any time $t$, $x \in \bbr^m$ and $\eta \in \bbr^m,$ we define 
$ {\cal V}_L (t,x, \eta) := \sum_{ \alpha : \| \alpha \| \le L- 1 } < (A_1)_{(\alpha)} (t,x), \eta >^2  $ and ${\cal V}_L (t, x) := \inf_{\eta : \| \eta \| = 1 } {\cal V}_L (t,x, \eta ) \wedge 1$. We are now ready to state our local weak H\"ormander condition at a given point $(t,y_0)$ where $t>0$ and $y_0\in E$:

{\bf (H4)} The H\"ormander condition is satisfied at $(t,y_0)$ if there exist $r\in ]0,t[\,$, $ R\in]0,1]\, $ and an integer $L$ such that $B_{5R} (y_0) \subset E$ and ${\cal V}_L (s,y)\geq c( y_0, R) > 0, \; \; \; \forall (s, y) \in [t-r,t]\times B_{3R} (y_0).$



\subsection{Local densities}\label{localdensities}

In this section we prove that ideas developed in \cite{stefano} can be extended to a time inhomogeneous SDE satisfying only a local weak H\"ormander condition. The non degeneracy assumption in \cite{stefano} is local ellipticity which fails to hold in our case. Let us remind the reader that an $\bbr^m $- valued random vector admits a density with respect to Lebesgue measure (or is absolutely continuous) on an open set $O \subset \bbr^m $ if for some function $ p \in L^1 (O)$ 
$$ E (f(X)) = \int f (x) p(x) dx ,$$
for any continuous function $f \in C_b ( \bbr^m ) $ satisfying $ supp (f) \subset O.$  We rely on the following classical criterion for smoothness of laws based on a Fourier transform method. 

\begin{prop}\label{fourier0}
 Let $\nu$ be a probability law on $\bbr^m $ and let $ \hat \nu (\xi)$ be its Fourier transform. If $\hat \nu $ is integrable, then $\nu$ is absolutely continuous and 
$$ p(y) = \frac{1}{(2 \pi)^m} \int_{\bbr^m} e^{ - i < \xi , y>} \hat \nu (\xi ) d \xi $$
is a continuous version of its density.  
\end{prop}

\begin{theo}\label{theo:main2}
Let $x \in E $ and $ ( X_t, t \geq 0 )$ be a strong solution of (\ref{eq:equationms}), starting from $x$, satisfying {\bf (H1)}-{\bf (H3)}.  Assume moreover that {\bf (H4)} is satisfied at $(t,y_0)$. Then the random variable $X_t $ admits a density $ p_{ 0, t} ( x, y ) $ on $ B_R ( y_0 ) $ which is continuous with respect to $ y \in B_R ( y_0 )  .$ 
\end{theo}
\begin{theo}\label{theo:lsc}Under the assumptions and notations of Theorem \ref{theo:main2}, for any fixed $ y \in B_R ( y_0) ,$ the map $x\in E \mapsto  p_{0, t } (x,y) $ is lower semi-continuous.
\end{theo} 

\noindent{\bf Proof of Theorem \ref{theo:main2}.} We have to use localization arguments. Let $\Phi $ be a function in $C_b^\infty ( \bbr^m ) $ such that $ 1_{ B_R (0)} \le \Phi \le 1_{ B_{2R}(0)} .$ Fix $x$ and $T$ such that $T\geq t$ and assume that $ E_x ( \Phi ( X_t - y_0 )) := m_0 > 0 .$ Then we can define a probability measure $\nu$ via 
\begin{equation}\label{eq:mu}
 \int f(y) \nu (dy) : = \frac{1}{m_0} E_x \left( f( X_t) \Phi ( X_t - y_0 )\right) 
\end{equation}
and its Fourier transform given by  
\begin{equation*}
\hat \nu (\xi ) = \frac{1}{m_0} E_x \left( e^{i < \xi , X_t>} \Phi ( X_t - y_0 ) \right). 
\end{equation*}
In order to show that $\hat \nu (\xi ) $ is integrable (cf. Proposition \ref{fourier0}) we use Malliavin calculus localized around $y_0.$ Let $\psi \in C^\infty_b (\bbr^m ) $ such that 
$$ \psi (y ) = \left\{ 
\begin{array}{ll} 
y & \mbox{ if } |y| \le 4 R \\
5 R \frac{ y}{|y|} & \mbox{ if } |y| \ge 5 R 
\end{array}
\right. $$
and $|\psi (y)| \le 5 R $ for all $ y .$ Let $ \bar b (t, y) = b (t, \psi ( y - y_0 ))$ and $\bar \sigma (y) = \sigma ( \psi ( y - y_0 ))$ be the localized coefficients of (\ref{eq:equationms}). Assumption {\bf (H2)} ensures that $\bar b$ and $\bar \sigma$ are $C_b^\infty-$extensions (w.r.t. $x$) of $ b_{| B_{4R} (y_0)} $ and $\sigma_{ | B_{4R} (y_0)} $ with $\bar b$ and its derivatives bounded on $[0,T]$. Let $ \XXX $ denote the unique strong solution of the equation 

\begin{equation}\label{eq:processgood}
\bar X_{i,s} = x_i + \int_0^s \bar b_i( u, \bar X_u) du + \int_0^s \bar \sigma_i (\bar X_u) d W_u ,\;  u \le T ,\; 1 \le i \le m ,
\end{equation}
with $\bar X_0=X_0=x$. If $ x \in B_{4R } ( y_0), $ the processes $ \XXX $ and $ X$ coincide up to the first exit time of $ B_{4R} (y_0)$. In the sequel we make use of the classical notation for flows $ \bar X_{s, t } (x) $, to denote the value at time $t$ of the solution starting at $x$ at time $s$ where $s\leq t$.

For a fixed $ \delta \in ] 0, t/2 \wedge r[,$ let $\tau_1 := \inf \{ s \geq t - \delta : X_s \in B_{3R} (y_0 ) \} $ and $\tau_2 := \inf \{ s \geq \tau_1 : X_s \notin B_{ 4 R }(y_0) \} .$ 
Then the set $\{ \Phi ( X_t - y_0 ) > 0 \} $ is equal to the union 
\begin{equation*}
\{ \Phi ( X_t - y_0 ) > 0 ; t - \delta = \tau_1 < t < \tau_2  \} 
\cup \left\{ \Phi ( X_t - y_0 ) > 0 ; \sup_{ 0 \le s \le \delta } |  \bar X_{\tau_1 , \tau_1 + s} ( X_{\tau_1 } ) - X_{ \tau_1} | \geq R \right\} .
\end{equation*}
Hence using the Markov property in $\tau_1$,
\begin{eqnarray*}\label{eq:fourier}
m_0 \hat \nu ( \xi ) &=&E_x \left( e^{ i < \xi, X_t>} \Phi ( X_t - y_0 )1_{\Phi ( X_t - y_0 ) > 0 ; \sup_{ 0 \le s \le \delta } | \XXX_{\tau_1, \tau_1 + s } ( X_{\tau_1}) - X_{ \tau_1} | \geq R }  \right)\\
&+& E_x \left( e^{ i < \xi, X_t>} \Phi ( X_t - y_0 )1_{ \Phi ( X_t - y_0 ) > 0 ; t- \delta = \tau_1 < t < \tau_2} \right) .
\end{eqnarray*}
The first term (call it A) on the right hand side of this equality is controlled by classical estimates. The important contribution comes from the second term (call it B). Indeed, using the following classical estimate 
\begin{equation}\label{eq:ub1}
E\left( \sup_{ u :s \le u \le t } | \bar X_{i,u} - \bar X_{i,s} |^q \right) \le C(T,q ,m,  b ,  \sigma ) (t-s)^{q/2} \mbox{ for all $ 0 \le s \le t \le T ,$}
\end{equation}
we can control A as follows. For all $ q > 0 , $
\begin{equation}\label{eq:control11}
P_x\left( \Phi ( X_t - y_0 ) > 0 ; \sup_{ 0 \le s \le \delta } | \XXX_{\tau_1, \tau_1 + s } ( X_{\tau_1}) - X_{ \tau_1} | \geq R \right) \le 
C(T,q,m , b, \sigma ) R^{-q} \delta^{q /2}  .
\end{equation}
The above estimation holds uniformly in $x.$ The constant $C(T,q ,m,  b ,  \sigma ) $ depends  on the supremum norms of $\bar b $ and $\bar \sigma,$ hence, by construction, on the supremum norms of $\sigma $ (resp. $b$) on $B_{5R} (y_0 )$ (resp. $B_{5R} (y_0 )\times [0,T]$).

The second term B can be controlled as follows. Thanks to the Markov property at time $ t - \delta,$ we first bound the modulus of B
\begin{equation}\label{eq:control2}
\Big| B \Big| \le 
 \sup_{ y \in B_{3R} (y_0)} | E \left( e^{ i < \xi, \XXX_{t-\delta, t} (y) > } \Phi ( \XXX_{t-\delta, t} (y)  - y_0 ) \right) | .
\end{equation}
Again this control holds uniformly in $x.$ In the right-hand side of (\ref{eq:control2}) we differentiate twice with respect to each space variable and then we apply the integration by parts formula of Malliavin calculus. Introducing the multi-index
$ \beta := ( 1,1, 2,2, \ldots , m,m ) $ and using the identity $ \partial_{x_k} e^{ i < \xi, x >}  = i \xi^k e^{ i < \xi , x > } , $ we obtain
\begin{eqnarray*}
| E \left( e^{ i < \xi , \XXX_{t-\delta, t} (y) >} \Phi ( \XXX_{t-\delta, t} (y)  - y_0) \right) | &\le& 
\Big({\prod_{l =1 }^m | \xi_l |^2 \Big)}^{-1} \Big| E \left( \partial_\beta e^{ i < \xi, \XXX_{t-\delta, t} (y)>} \Phi ( \XXX_{t-\delta, t} (y)- y_0 ) \right) \Big| \\&\le& 
\Big({\prod_{l =1 }^m | \xi_l |^2 \Big)}^{-1} \Big| E \left( | H_\beta ( \XXX_{t-\delta, t} (y) , \Phi ( \XXX_{t-\delta, t} (y) - y_0 ) ) \Big| \right ) .
\end{eqnarray*}
$H_\beta $ is the weight resulting from the integration by parts formula (see e.g.\cite{stefano}, Proposition 2.1). Remember that $L$ is the number of brackets needed in order to span $\bbr^m $ at any point of $B_{3R} (y_0 )$ (cf. $({\bf H4} ))$. We will show in the Appendix that the following classical result holds: there exists a constant $k_L $ such that 
\begin{equation}\label{eq:classical}
\|  H_\beta ( \XXX_{t-\delta, t} (y) , \Phi ( \XXX_{t-\delta, t} (y) - y_0 )\|_p \le C (r,  p, R, m)  \delta^{ - m k_L  } .
\end{equation}
We deduce from (\ref{eq:control11}) and (\ref{eq:classical}) that, for any $ q \geq 1$ and any $0 <  \delta < \frac{t}{2} \wedge r ,$  
$$
m_0 |\hat \nu (\xi )| \le C (T, r, R, q, m ) \;
 \left[ R^{- q } \delta^{ q/2} + \Big({\prod_{l=1 }^m | \xi_l |^2 \Big)}^{-1}\delta^{ - m k_L } \right] .
$$

The key point now is the freedom that still remains in the choice of $\delta$ and $q$. This is the main idea of balance given in \cite{stefano}: for a given $\xi$, choose $\delta$ and $q$ such that $ R^{- q } \delta^{ q/2} + \Big({\prod_{l=1 }^m | \xi_l |^2 \Big)}^{-1} \delta^{ - m k_L }$ tends to zero faster than $ \left( \prod_{l=1}^m  | \xi_l | \right)^{ - 3/2  }$ as $\| \xi \| \to \infty. $  More precisely here, the choice 
$$ \delta = t/2 \wedge r \wedge \| \xi \|^{-\frac{1}{2mk_L} }, \; q = 6 m  k_L  $$
where $ \| \xi \| := \prod_{l=1}^m  | \xi_l | ,$  leads to the integrability of $\hat \nu$ in $\xi $ for $\| \xi \| \to \infty$ since then
\begin{equation}\label{eq:ub5} 
m_0 |\hat \nu (\xi ) |\le C (T, r, R, q, m ) \; \| \xi \|^{ - \frac32}.
\end{equation}

Now we are able to conclude the proof of Theorem \ref{theo:main2}. From Proposition \ref{fourier0}, for any $ y \in B_R (y_0), $ 
\begin{equation}\label{eq:ouf}
p_{0, t } ( x, y ) = \frac{m_0}{(2 \pi)^m} \int_{\bbr^m} e^{ - i < \xi , y > } \hat \nu (\xi ) d \xi =
\frac{1}{(2 \pi)^m} \int_{\bbr^m } e^{ - i <\xi , y > } E_x ( e^{ i < \xi , X_t>} \Phi ( X_t - y_0 ) ) d \xi .
\end{equation}
We split the latter integration in two parts: over a finite region $I$ where $ \| \xi \| \le C $ and over its complement $I^c $. 
On $I,$ the integrand is bounded above by $ 1 $ (since $ \Phi \le 1_{ B_{2R} (0)} $), on $I^c $ we use the upper bound (\ref{eq:ub5}). This proves the continuity of $ p_{ 0, t} (x, y ) $ with respect to $y.$ This continuity is uniform in the starting point $x,$ since the upper bounds in (\ref{eq:control11}) and (\ref{eq:control2}) do not depend on $x.$\halmos

\noindent{\bf Proof of Theorem \ref{theo:lsc}.} In order to prove the lower semi-continuity of $p_{0, t} (x,y) $ in $x \in E ,$ for fixed $y \in B_R (y_0),$ we show that $x\mapsto p_{0, t} (x,y)$ is the limit of an increasing sequence of continuous functions $x \mapsto p_{0, t}^{(n)} (x,y) $. This will imply that it is lower semi-continuous. We also use localization arguments here but now we compare the diffusion $X$ with an approximation $X^{(n)}$ obtained by considering $X$ before the first exit time of some compact $K_n$ (cf. Assumption {\bf (H1})). It is then natural to use the flow property of $X^{(n)}$ which implies continuous dependence on the starting point. Note that the process $X$ itself may not enjoy the flow property.

We keep the notations introduced in the proof of Theorem \ref{theo:main2}. Given an integer $n$, let $b^{(n)} (t,x) $ and $\sigma^{(n)} (x)$ denote $C^\infty-$extensions (in $x$) of $ b(t, \cdot _{ | K_n} )$ and $\sigma_{| K_n} .$ Let $X^{(n)}$ be the associated diffusion process solving the localized SDE. The first exit time of $K_n$ is denoted by $T_n$ (cf {\bf (H1)}). Using that $T_n \to \infty$ and $X_t^{(n)} = X_t$ on $\{ T_n > t \} $ almost surely, we can write the following, for any $x \in K_n$ and any positive measurable function $f:$ 
\begin{eqnarray*}
\int f(y) \nu (dy) = \frac{1}{m_0} E_x \left( f(X_t)  \Phi ( X_t - y_0) \right) &=& \lim_{n \to \infty } \frac{1}{m_0} E_x \left( f( X_t)  \Phi (X_t - y_0 ) 1_{ \{ T_n > t \} } \right) \\
& \geq  &  \frac{1}{m_0} E_x \left( f( X_t ) \Phi (X_t - y_0 )  1_{ \{ T_n > t \}} \right) \\
&= & \frac{1}{m_0} E_x \left( f( X^{(n)}_t ) \Phi (X_t^{(n)} - y_0 )  1_{ \{ T_n > t \}} \right) .
\end{eqnarray*}

We approximate $1_{ \{ T_n > t \}} $ by some continuous functional on $ \Omega := C( \bbr_+ , \bbr^m ) .$  The set $\Omega $ is endowed with the topology of uniform convergence on compacts. $\mathbb{P}_{0,x}^{(n)}$ denotes the law of $X^{(n)} $ on $ (\Omega , {\cal B} (\Omega ) ),$ starting from $x$ at time $0.$ We know that the family $\{ \mathbb{P}_{0,x}^{(n)}, x \in \bbr^m \}$ is Feller, i.e. if $x_k \to x, $ then $\mathbb{P}_{0, x_k}^{(n)} \to \mathbb{P}_{0,x}^{(n)} $ weakly as $k \to \infty .$  Thanks to this property, $ m_0^{ - 1} E_x \left( f( X^{(n)}_t) \Phi (X_t^{(n)} - y_0 )   \right)$ is continuous w.r.t. $x$. Let $M_t^n = \max_{ s \le t } X_s^{(n)} $ and $m_t^n = \min_{s \le t} X_s^{(n)} $ be the (coordinate-wise) maximum and minimum processes associated to $X^{(n)}.$  Due to the structure of the compacts $K_n$ (see assumption {\bf (H1)}), we can construct  
$C^\infty-$functions $ \varphi^n, \Phi^n $ such that $ 1_{ [a_{n-1}, \infty [ } \le \varphi^n \le 1_{ [a_n,  \infty[ } $ and 
$ 1_{ ] - \infty , b_{n-1}] } \le \Phi^n \le 1_{ ] - \infty , b_n] } $ (these inequalities have to be understood coordinate-wise). Then, since $X_t$ equals $X_t^{(n)} $ up to time $T_n,$ 
$$ \{ T_{n-1} > t \} = \{ a_{n -1} \le m_t^n \le M_t^n \le b_{ n - 1} \} \subset \{ \varphi^n( m_t^n ) > 0 , \Phi^n (M_t^n ) > 0  \} \subset \{ T_n > t \} . $$ 
So for any $ f \geq 0 , $ 
\begin{eqnarray*}
E_x \left(f( X^{(n)}_t)  \Phi (X_t^{(n)} - y_0 ) 1_{ \{ T_n > t \} } \right) 
& \geq & E_x \left( f( X^{(n)}_t) \Phi (X_t^{(n)} - y_0 )  \Phi^n ( M_t^n) \varphi^n (m_t^n) \right).
\end{eqnarray*}
Define now a sub-probability measure $\nu_n $ by 
\begin{equation}\label{eq:mun}
 \int f(y) \nu_n ( dy ) := \frac{1}{m_0}  E_x \left( f( X^{(n)}_t ) \Phi (X_t^{(n)} - y_0 )  \Phi^n ( M_t^n) \varphi^n (m_t^n ) \right). 
\end{equation}
The new functional 
$$ \Phi (X_t^{(n)} - y_0 )  \Phi^n ( M_t^n) \varphi^n (m_t^n ) $$
satisfies the same hypotheses as the former $ \Phi (X_t^{(n)} - y_0 ) . $
For any $f \geq 0 , $ we have that 
$$ \int f(y) \nu_n (dy) \le \int f(y) \nu_{n+1} (dy) \uparrow \int f(y) \nu (dy ) \mbox{ as } n \to \infty .$$ 
If we can show that $ \nu_n $ possesses a density, that we shall denote by $m_0^{ - 1 }  p_{0, t}^{(n)} (x, y ) ,$ the following inequalities will hold true  
\begin{equation}\label{eq:increase}
 p_{0, t}^{(n)} (x, y ) \le p_{0, t}^{(n+1)} (x, y ) \le p_{0, t} (x, y ) \quad \mbox{ for all $ n \geq 1 ,$}  
\end{equation}
for any fixed $ x ,$ $\lambda (dy) -\mbox{almost surely.}$
So in a next step we show that indeed $\nu_n $ possesses a density. In order to indicate explicitly the dependence on the starting point $x$, we introduce
$$\gamma_n ( x, \xi) := \hat \nu_n ( \xi ) =   \frac{1}{m_0}  E_x \left( e^{i < \xi, X^{(n)}_t>} \Phi (X_t^{(n)} - y_0 )  \Phi^n ( M_t^n) \varphi^n (m_t^n ) \right) $$
for the Fourier transform of $\nu_n$ and we apply the argument in the proof of Theorem \ref{theo:main2} to $ \gamma_n ( \cdot, \xi ).$ The upper bounds (\ref{eq:control11}), (\ref{eq:control2}) and (\ref{eq:ub5}) hold also for $m_0 \gamma_n (x, \xi).$ Moreover, they hold uniformly in $x .$ This implies first that $\xi \to \gamma_n ( x, \xi ) $ is integrable. Hence, the density of $ m_0 \nu_n $ exists and is given by
\begin{equation}\label{eq:pn}
p_{0, t}^{(n)} (x, y ) =  \frac{m_0}{(2\pi)^m} \int_{\bbr^m } e^{ - i < \xi , y >} \gamma_n ( x, \xi ) d \xi .
\end{equation}
From the common upper bounds (\ref{eq:control11}), (\ref{eq:control2}) and (\ref{eq:ub5}) and the fact that $ \gamma_n ( x, \xi ) \to \hat \nu ( \xi )$ as $ n \to \infty$ we obtain that
$$p_{0, t}^{(n)} (x, y ) \to p_{0, t}  (x, y ) ,$$
Taking into account  (\ref{eq:increase}), this implies 
$$ p_{0, t} (x,y) = \lim_n \uparrow p_{0, t}^{(n)} (x,y) .$$

It remains to show that for any $y \in B_R ( y_0),$ the map $x \mapsto p_{0, t}^{(n)} (x,y) $ is continuous. This follows from the continuity of $\gamma_n ( x, \xi) $ in $x $ and the fact that the upper bounds (\ref{eq:control11}), (\ref{eq:control2}) and (\ref{eq:ub5}) hold uniformly in $x,$ by dominated convergence. The continuity of $\gamma_n ( x, \xi) $ in $x$ follows from the Feller property of $\mathbb{P}_{0,x}^{(n)} $ and the fact that all operations appearing in $ \gamma_n (x, \xi) $ are continuous on $\Omega.$  \halmos

\section{Densities for SDE with internal variables and random input}\label{ivri}

In this section we consider the stochastic systems presented in Section \ref{sec:internal} with internal variables and random input
 \begin{eqnarray}\label{internal}
dX_{1,t}&=& F (X_{l,t}, \ 1\leq l\leq m-1)dt +dX_{m,t}\, ,\\
dX_{i,t}&=&[-a_i(X_{1,t})X_{i,t}+b_i(X_{1,t})]dt, \, \, \, i = 2, \cdot, \cdot, \cdot , m-1\, ,\nonumber\\
dX_{m,t} &=& b_m(t,X_{m,t})dt \;+\sigma (X_{m,t})dW_t\, ,\nonumber 
 \end{eqnarray}
under the assumptions {\bf (H1)}-{\bf (H3)}. Let us define $E_m:= \bbr\times[0,1]^{m-2}\times U$ where $U$ is the interval where $X_m$ evolves (cf. assumption {\bf (H3)}). We first express a sufficient condition for the weak H\"ormander condition {\bf (H4)} to be satisfied and hence the existence of densities to hold locally, based on the previous section. Then we address the question of positivity of these densities. 

\subsection{Weak H\"ormander condition for SDE with internal variables and random input}\label{weakHforivri}

\begin{defin}\label{determinant}
Let $k$ be an integer and denote by $\partial^{(k)} _{x_1}$ the partial derivative w.r.t. $x_1$ of order $k$. For any $x\in E_m$ define the column vector $J(x)\in \bbr^{m-1}$ by $J_1(x):=F(x_1,x_2, \cdot,\cdot,\cdot,x_{m-1})$, $J_i(x):=-a_i(x_1)x_i+b_i(x_1)$, $\ 2\leq i\leq m-1$. Then ${\bf D}(x)$ is the determinant of size $m-1$, whose columns are the $\partial^{(k)} _{x_1}J(x)$ with $1\leq k\leq m-1$.
\end{defin}

\begin{theo}\label{theo:internal_hoer}
The weak H\"ormander condition {\bf (H4)} holds at any point $x=(x_i, \ 1\leq i\leq m)\in E_m$ such that ${\bf D}(x) \ne 0$. 
\end{theo}
It is important to note that ${\bf D}(x)$ actually depends only on the $m-1$ first components of $x$. In particular if the $m-1$ first components of two points $x$ and $x'$ coincide, then ${\bf D}(x)={\bf D}(x')$. This remark will be important in the sequel. Moreover the condition in Theorem \ref{theo:internal_hoer} implies a version of {\bf (H4)} uniform w.r.t. time on every compact interval $[0,T]$.

In the proof of Theorem \ref{theo:internal_hoer} we will rely on the following Proposition about computation of Lie brackets.
\begin{prop}\label{calculcrochetsdelie;internalvariables2}
We identify $x_0$ with $t$. Let $e\in \bbr^{m+1}$ be defined by $e_i=0$ for $i=0$ and $2\leq i\leq m-1$  and $e_1=e_m=1$. Let $\Xi(t,x):=\varphi(x_m)e$ for some function $\varphi$. Let $Y(t,x)$ such that $Y=\psi(x_m)e+\rho(x_m)\sum_{i=1}^{m-1}y_i\frac{\partial}{\partial x_i}$ where the functions $y_i$ satisfy $\partial^{(1)}_{x_m}y_i\equiv 0$ for $1\leq i\leq m-1$. Then for all $(t,x)$, the Lie bracket $[\Xi,Y]$ satisfies
\begin{equation}
[\Xi,Y](t,x)=(\varphi \psi'-\varphi' \psi)(x_m) e+\varphi(x_m)\rho'(x_m)\sum_{i=1}^{m-1}y_i\frac{\partial}{\partial x_i}+\varphi(x_m)\rho(x_m)\sum_{i=1}^{m-1}\partial_{x_1}^{(1)}y_i\frac{\partial}{\partial x_i}.
\end{equation}
\end{prop}

\noindent {\bf Proof of Proposition \ref{calculcrochetsdelie;internalvariables2}} By assumption, 
\begin{eqnarray*}
[\Xi, Y]&=&[\varphi(x_m)e,\psi(x_m)e+\rho(x_m)\sum_{i=1}^{m-1}y_i\frac{\partial}{\partial x_i}]\\
&=&[\varphi(x_m)e,\psi(x_m)e]+[\varphi(x_m)e,\rho(x_m)\sum_{i=1}^{m-1}y_i\frac{\partial}{\partial x_i}]\\
&=&(\varphi \psi'-\varphi' \psi)(x_m) e+[\varphi(x_m)e,\rho(x_m)\sum_{i=1}^{m-1}y_i\frac{\partial}{\partial x_i}]\\
&=&(\varphi \psi'-\varphi' \psi)(x_m) e+\varphi(x_m)\sum_{i=1}^{m-1}[\rho'(x_m)y_i+\rho(x_m)\partial_{x_1}^{(1)}y_i]\frac{\partial}{\partial x_i}.\quad \quad\quad \quad\halmos
\end{eqnarray*}

\noindent {\bf Proof of Theorem \ref{theo:internal_hoer}.} We keep the notations of Section \ref{sec:hoerm} as well as those of Proposition \ref{calculcrochetsdelie;internalvariables2}. Let $A_1=\sigma(x_m)(\frac{\partial}{\partial x_1}+\frac{\partial}{\partial x_m})=\sigma(x_m)e$ and $A_0=\sum_{i=0}^m A_{0,i}\frac{\partial}{\partial x_i}$ with $A_{0,0}\equiv 1$. We compute the Lie brackets defined recursively by $L_1:= [A_1,A_0]$ and $L_{k+1}=[A_1,L_k]$. Let us start with $L_1$. By definition $L_1=A_1 A_0-A_0 A_1$. This implies
\begin{equation*}
L_1=\sigma(x_m)\sum_{i=1}^m\frac{\partial A_{0,i}}{\partial x_1}\frac{\partial}{\partial x_i}+\sigma(x_m)\frac{\partial A_{0,1}}{\partial x_m}\frac{\partial}{\partial x_1}+\sigma(x_m)\frac{\partial A_{0,m}}{\partial x_m}\frac{\partial}{\partial x_m}
-\sigma'(x_m)A_{0,m}(\frac{\partial}{\partial x_1}+\frac{\partial}{\partial x_m}).
\end{equation*}
Now in the system of interest, $\frac{\partial A_{0,m}}{\partial x_1}\equiv 0, $ $\frac{\partial A_{0,1}}{\partial x_m}\equiv \frac{\partial A_{0,m}}{\partial x_m}$ and the $\frac{\partial A_{0,i}}{\partial x_1}, 2\leq i\leq m-1, $ do not depend on $x_m$. Therefore $L_1$ satisfies the assumptions of Proposition \ref{calculcrochetsdelie;internalvariables2} which applies with $\Xi=A_1$ and $Y=L_1$. Moreover $\varphi(x_m)\equiv\rho(x_m)=\sigma(x_m), $ $y_i=\frac{\partial A_{0,i}}{\partial x_1}$ for $1\leq i\leq m-1,\ $ and $\quad  \psi(x_m)=\sigma(x_m)\frac{\partial A_{0,1}}{\partial x_m}-\sigma'(x_m)A_{0,m}$. Hence
\begin{equation}\label{secondbracket}
L_2=[A_1,L_1]=\Psi_1(x_m)e+\Psi_2(x_m)\sum_{i=1}^{m-1}\partial_{x_1}^{(1)} A_{0,i}\frac{\partial}{\partial x_i}+\Psi_3(x_m)\sum_{i=1}^{m-1}\partial_{x_1}^{(2)} A_{0,i}\frac{\partial}{\partial x_i}\, ,
\end{equation}
where the functions $\Psi_l$ are expressed using $\varphi, \rho, \psi$ and their derivatives. Again, identity (\ref{secondbracket}) coupled with Proposition \ref{calculcrochetsdelie;internalvariables2} enables us to work by iteration. We obtain the expression of $L_k$ for any $k\geq 1$:
\begin{equation}\label{kbracket}
L_k=\Phi_1(x_m)e+\sum_{i=1}^{m-1}\sum_{l=1}^{k-1}\Phi_l(x_m)\partial_{x_1}^{(l)} A_{0,i}\frac{\partial}{\partial x_i}+\Psi_{k+1}(x_m)\sum_{i=1}^{m-1}\partial_{x_1}^{(k)} A_{0,i}\frac{\partial}{\partial x_i}.
\end{equation}
The statement of Theorem \ref{theo:internal_hoer} follows immediately from (\ref{kbracket}) thanks to the multilinear property of determinants.
\halmos

\begin{defin}
Define the set ${\cal H} := \{ (x_i, \ 1\leq i\leq m)\in \bbr^{m-1}\times U; \quad{\bf D}(x) \ne 0\}.$
\end{defin}


The set ${\cal H}$ is an open subset of $E_m$ by continuity of ${\bf D}$ on $\bbr^{m-1}\times U$.
 
\begin{theo}\label{localdensityinternal}
Assume that there exists $y_0\in {\cal H}$ and $R>0$ such that $B_{3R}(y_0)\subset {\cal H}$. Then for any $x\in E_m$ and $t>0$, the random variable $X_t $ admits a density $ p_{ 0, t} ( x, \cdot ) $ on $ B_R (y_0 ) $. Moreover the map $ y \in B_R (y_0)\mapsto p_{ 0, t} ( x, y )$ is continuous  and for any fixed $ y \in B_R (y_0) ,$ the map $x \in E_m \mapsto  p_{0, t} (x,y) $ is lower semi-continuous.
\end{theo}

\noindent {\bf Proof.} The statement is a direct consequence of Theorems \ref{theo:main2} and \ref{theo:lsc} of Section \ref{localdensities}. \halmos

The difficulty in practice is to obtain more information on the set ${\cal H}$,  in particular to know whether it coincides with $E_m$ for a given system. At least one would like to be able to specify open regions included in ${\cal H}$. One can hope to achieve this goal only numerically unless the coefficients of the system are very simple. In Section \ref{sec:main} we provide details for a stochastic Hodgkin-Huxley model.

The definition of ${\cal H}$ obeys the intuition that in (\ref{internal}) the noise is most rapidly transported through the first and the last variables (the only ones carrying Brownian noise). In some sense, when proceeding in this way, we are exploring directions of the space where the diffusion moves at maximal possible speed. Accordingly, except the first Lie bracket $L_1= [A_0,A_1] $ which involves the drift $A_0$, we always use the diffusion coefficient $A_1$ in order to compute the brackets of higher order. Developing the solution $X$ of (\ref{internal}) for small time steps $\delta$ into iterated Ito integrals shows that the speed of the diffusion is of order $ \delta^{\frac12}$ in the direction of $A_1$, of order $ \delta^{ 1 + \frac12}$ in the direction of $L_1$ (the $+1$ comes from the drift $A_0$). For the subsequent Lie brackets each time we use $A_1$ we add a factor $\frac{1}{2}$ so that the speed of the diffusion is of order $\delta^2$ in the direction of $L_2$, of order $\delta^{1 + 3 \times \frac12}$ in the direction of $L_3$  and so on. We refer the reader to \cite{NSW}, in particular identity (12).
Hence it is important to remember that belonging to ${\cal H}$ is only a sufficient condition for the weak H\"ormander condition to hold. It may hold also at points outside ${\cal H}$ in which case the system suffers a slow down in the sense just explained.

\subsection{Positivity of densities}\label{positivedensities}
Once we have proved that densities exist, even if only locally, we look for regions where they are positive. For this purpose we combine control arguments and the support theorem. We keep the notation $E_m$ from the previous section. We start by proving two accessibility results for (\ref{internal}) which hold without any assumption on the existence of densities. These two results are different in nature since the first one relies on stability properties of the underlying deterministic system (\ref{subinternal}) while the second one does not. For the proof of the first result (Proposition \ref{theo:accessible}) we refer the reader to \cite{BLBMZ-13} for similar ideas in the framework of Piecewise Deterministic Markov Processes.

Below we denote by $P_{0, t } ( x , \cdot )$ the law of $X_t$ which solves (\ref{internal}) with initial value $x$. We make use of the notations introduced in Proposition \ref{globallyasymptoticallystable} and the subsequent remark: $y_{i,\infty}(z_1)=\frac{b_i(z_1)}{a_i(z_1)},\, \,i = 2, \cdot, \cdot, \cdot , m-1,$ for any given $z_1$ where the functions $a_i$, $b_i$ appear in (\ref{subinternal}) that we recall here for the reader's convenience:
\begin{eqnarray*}
dz_{1,t}&=& F (z_{l,t}, \ 1\leq l\leq m-1) dt \, , \\
dz_{i,t}&=&[-a_i(z_{1,t})z_{i,t}+b_i(z_{1,t})]dt, \, \, \, i = 2, \cdot, \cdot, \cdot , m-1.
 \end{eqnarray*}

\begin{prop}\label{theo:accessible}
Let $z_1$ be an arbitrary real number. Define $\overline z\in \bbr^{m-1}$ by $\overline z:= (z_1, y_{i,\infty}(z_1)\, ,\,i = 1, \cdot, \cdot, \cdot , m-1)$. For all $x\in E_m$ and any neighborhood ${\cal N}$ of $ \overline z$ in $\bbr\times(0,1)^{m-2}$ there exists $t_0 $ such that  
\begin{equation}
 \forall t \geq t_0,\quad P_{0, t } ( x , {\cal N} \times U ) > 0 .
\end{equation}

This statement holds in the particular case where $F(\overline z)=0$. In this case the point $\overline z$ is an equilibrium point of (\ref{subinternal}). 
\end{prop} 
\noindent{\bf Proof of Proposition \ref{theo:accessible}.}\label{sec:supp} As in the proof of Theorem \ref{theo:lsc} we write 
$ \Omega$ for $C ( [ 0, \infty [ , \bbr^m )  $
and endow it with its canonical filtration $ ( {\cal F}_t)_{t \geq 0 } .$ Recall that $ \mathbb{P}_{0,x}$ is the law of $(X_{  u }, u \geq 0 ), $ starting from $x$ at time $0.$ We first localize the system by a sequence of compacts $(K_n)$ and let $T_n = \inf \{ t : X_t \in K_n^c \} $ be the exit time of $K_n$ (cf. assumption {\bf (H2)}). For a fixed $n,$ let $ b^{(n)} (t, x) $ and $\sigma^{(n)} (x)$ be $C_b^\infty -$extensions in $x$ of $ b(t, \cdot_{| K_n }) $ and $\sigma_{| K_n}$ (here we denote by $b$ and $\sigma$ the coefficients of (\ref{internal}) for short). Let $X^{(n)}$ be the associated diffusion process. For any integer $n\geq 1$ and starting point $x,$ we write $\mathbb{P}_{0, x}^{(n )}$ for the law of $(X^{(n)}_{  u } , u \geq0  )$ on $\Omega$ satisfying $X_0=x$. We wish to find lower bounds for quantities of the form $\mathbb{P}_{0,  x } (  B ) $ where $B = \{ f \in \Omega : f(t) \in {\cal N}\times U \}  \in {\cal F}_t$, for given $t>0$. We start with the following inequality which holds for any $t > 0 $ and $n$:
\begin{equation}\label{eq:tobelb}
\mathbb{P}_{0,  x } (  B ) \geq \mathbb{P}_{ 0, x } ( \{ f  \in B ; T_n > t \}  ) 
 = \mathbb{P}^{(n)}_{ 0, x } (  \{ f \in B ; T_n >  t \} )  .
\end{equation}
In the sequel we show that for some integer $n_0$ and any fixed $x \in K_{n_0}$, the quantity $\mathbb{P}^{(n)}_{ 0, x } (  \{ f \in B ; T_n >  t \} )$ is indeed positive provided that $n$ is sufficiently large. We are therefore interested in the support of $\mathbb{P}^{(n)}_{ 0, x }$. Fix $t$ and let $ {\cal C} := \{ {\tt h} : [ 0, t ] \to \bbr : {\tt h}(s) = \int_0^s \dot {\tt h} (u) du , \forall s \le t , \int_0^t \dot {\tt h}^2 (u) du < \infty \} $ be the Cameron-Martin space. 
Given ${\tt h} \in {\cal C} ,$ consider $ X({\tt h})\in \bbr^m$ the solution of the differential equation 
\begin{equation}\label{eqcontrol1}
X({\tt h})_s = x + \int_0^s \sigma^{(n)} ( X({\tt h})_u) \dot {\tt h} (u) du + \int_0^s  {\tilde b}^{(n)} ( u, X({\tt h})_u) du ,\quad s\le t.
\end{equation}
If (\ref{eqcontrol1}) were time homogeneous, the support theorem would imply that the support of $\mathbb{P}^{(n)}_{0, x}$ in restriction to $ {\cal F}_t$ is the closure of the set $ \{ X({\tt h}) : {\tt h} \in {\cal C} \} $ with respect to the uniform norm on $ [ 0, t]$ (see e.g.  \cite{Millet-Sanz} Theorem 3.5 or \cite{BenArous} Theorem 4). To conclude in our situation as well, it is enough to replace the $m-$dimensional process $X^{(n)}$ by the $(m+1)-$dimensional process $(t, X^{(n)} _t) $ which is time-homogenous. 
In order to proceed further we construct a control ${\tt h}$ so that $X({\tt h})$ remains in $K_n$ during $ [ 0, t ]$ provided that $n$ is sufficiently large.\\ 
We start by exploiting stability properties of the underlying deterministic system (\ref{subinternal}). Let $ \gamma: \bbr\mapsto\bbr$ a smooth function satisfying $\gamma(t):=z _1$ for all $t\geq 1$. Consider $y_t\in\bbr^{m-2}$ solving
\begin{equation*}
dy_{i,s}=[-a_i(\gamma(s))y_{i,s}+b_i(\gamma(s))]ds, \, \, \, i = 2, \cdot, \cdot, \cdot , m-1.
\end{equation*}
Then for all $t>1$, 
\begin{eqnarray*}
y _{i,t} &=& y_{i,0} e^{ - \int_0^t  a_i ( \gamma(s)) ds } + \int_0^t b_i ( \gamma(u)) e^{ - \int_u^t a_i( \gamma(r)) dr } du\\
&=&y_{i,1}e^{ - a_i ( z_1)(t-1)}+y_{i,\infty}(z_1)(1-e^{ - a_i ( z_1)(t-1)})
\end{eqnarray*}
where $y_{i,\infty}(z_1^*)=\frac{b_i ( z_1)}{a_i ( z_1)}.$ This formula expresses the fact that on $[1,+\infty[$, the coefficients $a_i(\gamma(s))$  (resp. $b_i(\gamma(s))$) are constant equal to $a_i ( z_1)$ (resp. $b_i ( z_1)$). Hence for any  $\varepsilon>0$, like in Proposition \ref{globallyasymptoticallystable}, there exists $t_0>1$ such that $|y_{i,t}-y_{i,\infty}(z_1)|<\varepsilon$ for all $t\geq t_0$ and all $2\leq i\leq m-2$. Now take $\varepsilon$ so small that $ B_\varepsilon ( \overline z) \subset {\cal N} .$ Then for all $t\geq t_0>1, $ the vector $(\gamma(t), y _{i,t}\, 2\leq i\leq m-2)$ belongs to $ B_\varepsilon ( \overline z ) $ (remember that for $t>1$, $\gamma(t)$ is fixed at $z_1$). \\ 
Fix an integer $n_0$ and $x$ in $K_{n_0}$. We are now able to construct a control $h\in{\cal C}$ such that the solution of (\ref{eqcontrol1}) remains in $K_n$ for all large enough times and all integer $n$ large enough. Choose a function $\gamma$ as above satisfying moreover $ \gamma (0 ) = x_1$, $ \gamma (1 ) =z _1$. Define $(Z_t)_{t\geq 0}\in \bbr^m$, the deterministic path starting from $x$ such that 
\begin{eqnarray}\label{construction}
Z_{1,s}&=&\gamma(s)\, ,  \\
dZ_{i,s}&=&[-a_i(Z_{1,s})Z_{i,s}+b_i(Z_{1,s})]ds, \, \, \, i = 2, \cdot, \cdot, \cdot , m-1\, ,\nonumber\\
Z_{m,s}&=& x_m-x_1+\gamma(s)-\int_0^s F(Z_u)du.\nonumber
 \end{eqnarray}
Then fix $t \geq t_0 $ and consider  a function ${\tt h}$ defined by
\begin{equation}\label{hdot}
 \dot {\tt h} (s) := \frac{\dot\gamma(s)- F(Z_s )-b_m(s,Z_{m,s})+\frac{1}{2}\sigma(Z_{m,s}) \sigma'(Z_{m,s})} {\sigma(Z_{m,s})}.
\end{equation}
Note that $(Z_s,\, s\in[0,t])$ is bounded and therefore remains in the compact $K_n$ for all $n$ large enough. In particular, $Z_{m,s}$ lies in a compact interval included in $U$ for all $s\leq t$. Then under assumption {\bf (H3)}, the expression (\ref{hdot}) is well-defined. This assumption also provides that $\dot{\tt h} \in L^2 ( [ 0, t ] ) ,$ hence ${\tt h} \in {\cal C} .$ Hence, with such a choice of ${\tt h},$ the solution $X({\tt h})$ of equation (\ref{eqcontrol1}) coincides with the solution $Z$ of system (\ref{construction}). As explained previously, we can choose $n$ such that $(Z_s,\, s\in[0,t])$ remains in $K_n$. \\
Consider now, for $\delta>0$, the tubular neighborhood $T_\delta$ of $(Z_s,\, s\in[0,t])$ in $\Omega$ of size $\delta$, namely the set  $\{ f \in  \Omega  : \sup_{ s \le t } | f(s) - Z_s | <  \delta \} .$ By the support theorem $ \mathbb{P}^{(n)}_{0, x} (T_\delta) > 0 .$ Remember that we have chosen $\epsilon$ and $t_0$ in order to satisfy $T_\delta \subset   \{ f \in \Omega : f(t) \in B_\varepsilon (\overline z) \times U \}    $ as well as $B_\varepsilon (\overline z)\subset{\cal N}$. Choosing $ \delta \le \varepsilon/2  $ such that $ T_\delta\subset \{ f \in \Omega : T_n (f) > t \} $, we conclude as announced that 
\begin{equation*}
 P_{0,t}(x,{\cal N}\times U)\geq  P_x ( X_t \in B_\varepsilon (\overline z) \times U  ) \geq \mathbb{P}^{(n)}_{0, x } (T_\delta) > 0 . 
 \end{equation*}
\halmos

The following statement is a consequence of Proposition \ref{theo:accessible} when one assumes in addition that densities exist. We keep the notations and assumptions of this proposition.
\begin{prop}\label{cor:3.1} 
Let $z_1$ be an arbitrary real number. Assume there exists $u\in U$ such that ${\bf D} (\overline z,u)\neq 0$. Define $y:=(\overline z,u)$. Then, for all $x\in E_m$,  there exist $\delta>0$ such that   
$$\forall \, t\geq t_0\, , \, \, \, \inf_{ u_1\in K_1 } \inf_{u_2\in K_2}  p_{ 0, t} (u_1, u_2 )> 0 $$
where $ K_1$ (resp. $K_2$) stands for the closure of $B_\delta (x)$ (resp. $B_\delta (y)$). 
\end{prop}

Before proving this statement let us notice that if ${\bf D} (\overline z,u)\neq 0$ for some $u\in U$ then ${\bf D} (\overline z,\tilde u)\neq 0$ for every $\tilde u\in U$ (see Theorem \ref{theo:internal_hoer}). 

\noindent{\bf Proof of Proposition \ref{cor:3.1}.}\label{sec:supp3} From Theorem  \ref{theo:main2} we see that the weak H\"ormander condition at $y$ implies the existence of densities locally around this point for all initial point $x$. The regularity of densities proved in Theorems \ref{theo:main2} and \ref{theo:lsc} combined with Proposition \ref{theo:accessible} yields the statement.\halmos

We now show that, during any arbitrary long period, with positive probability, the stochastic system (\ref{internal}) is able to reproduce the behavior of $(z_t,I(t))\in \bbr^m$ where $I$ is an arbitrary smooth input and $z_t$ solves (\ref{subinternal}) submitted to this input, namely
\begin{eqnarray}\label{subinternalwithI}
dz_{1,t}&=& [F (z_{l,t}, \ 1\leq l\leq m-1) +I(t)]dt  \\
dz_{i,t}&=&[-a_i(z_{1,t})z_{i,t}+b_i(z_{1,t})]dt, \, \, \, i = 2, \cdot, \cdot, \cdot , m-1.\nonumber
 \end{eqnarray}

Remember that $\Omega=C([0,\infty[;\bbr^m)$, $ \, B_{\delta}( x)$ denotes the open ball of radius $\delta$ centered at $x$.
\begin{prop}\label{theo:posinternal}
Fix $x\in E_m$ and $t>0$. Let $I$ be a smooth deterministic input such that $x_m+\int_0^s I(r)dr\in U$ for all $s\leq t$. Define $\mathbb{X}_s^x:=(\mathbb{Y}_s^{\tilde x}, \, x_m+\int_0^s I(r)dr, \, s\leq t)$ where $\mathbb{Y}^{\tilde x}$ is the deterministic path solution of (\ref{subinternalwithI}) starting from $\tilde x:=(x_i, 1\leq i\leq m-1)$. We denote by $\mathbb{P}_{0, x}$ the law of the solution of (\ref{internal}) starting at $x$. Then for any $\varepsilon > 0$ 
\begin{equation*}
\mathbb{P}_{ 0,x} \left( \left\{ f \in \Omega : \sup_{ s \le t } | f(s) - \mathbb{X}^x_s | \le \varepsilon \right\} \right) > 0 
\end{equation*}
and moreover there exists $\delta > 0$ such that for all $ x'' \in B_\delta ( x)$ 
\begin{equation*}
 \mathbb{P}_{ 0,x''} \left( \left\{ f \in \Omega : \sup_{ s \le t } | f(s) - \mathbb{X}^{x }_s | \le \varepsilon \right\} \right) > 0 .
 \end{equation*}
\end{prop}

\noindent{\bf Proof of Proposition \ref{theo:posinternal}.}\label{sec:supp2} We keep the notations introduced in the proof of Proposition \ref{theo:accessible}. In the course of this proof we have shown that the support theorem applies to inhomogeneous diffusions like the one obtained after localizing (\ref{internal}). Moreover we still hope to reach the positivity we are looking for through inequalities (\ref{eq:tobelb}) and paths solving (\ref{eqcontrol1}) for $h\in {\cal C}$, that remain in $K_n$ during $ [0, t]$ provided that $n$ is sufficiently large.  So the system we work with is the localized one. Consider $I$ to be a deterministic input such that $x_m+\int_0^s I(r)dr\in U$ for all $s\leq t$. Define $\chi_{m,s}:=x_m+\int_0^s I(r)dr$ for all $s\leq t$ and 
\begin{equation}\label{hdot2}
 \dot {\tt h} (s) := \frac{I(s)-b_m(s,\chi_{m,s})+\frac{1}{2}\sigma(\chi_{m,s})\sigma'(\chi_{m,s})}{\sigma(\chi_{m,s})}.
\end{equation}
By definition  $(\chi_{m,s}, s\leq t)$ lies in a compact interval included in $U$. Then, the expression (\ref{hdot2}) is well-defined by assumption {\bf(H3)}. This assumption also provides that $\dot{\tt h} \in L^2 ( [ 0, t] )$ hence ${\tt h} \in {\cal C}$. Moreover, with such a choice of ${\tt h},$ the controlled path $X({\tt h})$, solution of (\ref{eqcontrol1}), coincides with $(\mathbb{Y}_s^{\tilde x}, \, \chi_{m,s}, \, s\leq t)$ where $\mathbb{Y}^{\tilde x}$ is the deterministic path solution of (\ref{subinternalwithI}) starting from $\tilde x=(x_i, 1\leq i\leq m-1)$. We can choose $n$ large enough such that $(\mathbb{Y}_s^{\tilde x}, \, \chi_{m,s},\, s\in[0,t])$ remains in $K_n$. We write $\mathbb{X}_s^x$ for $(\mathbb{Y}_s^{\tilde x}, \, \chi_{m,s})$. 
For $\delta>0$, consider the tubular neighborhood $T_\delta$ of $\mathbb{X}^x$ on $[0,t]$ namely the set  $\{ f \in  \Omega  : \sup_{ s \le t } | f(s) - \mathbb{X}^x_s | <  \delta \} .$ By the support theorem $ \mathbb{P}^{(n)}_{0, x} (T_\delta) > 0.$ Choose now $ \delta$ such that $ T_\delta\subset \{ f \in \Omega : T_n (f) > t \} $. Taking $T_\delta$ as the set $B$ in (\ref{eq:tobelb}) yields the first statement of Proposition \ref{theo:posinternal}. The second one follows from the Feller property of $\mathbb{P}^{(n)}_{0, x }$ which enables us to extend the first statement to a small ball around $x$.
\halmos

We close this section with the following consequence of Proposition \ref{theo:posinternal} from which we borrow the notations.

\begin{prop}\label{cor:3} 
Let $x_m$ be an element of $U$ and $x:=(z^*,x_m)$ and $x':=(z^*, \, x_m+\int_0^t I(r)dr)$ which both belong to $E_m$ and assume moreover that ${\bf D} (x)\neq 0$. There exists $\delta>0$ such that  
$$ \inf_{ y\in K } \inf_{y' \in K'}  p_{ 0, t} (y, y') > 0 $$
where $ K$ (resp. $K'$) stands for the closure of $B_\delta (x)$ (resp. $B_\delta (x')$). 
\end{prop}

\noindent{\bf Proof of Proposition \ref{cor:3}.} Note that, from Theorem \ref{theo:internal_hoer}, the assumption ${\bf D}( x)\neq 0$ ensures that the local weak H\"ormander condition holds at both points $x$ and $x'$. \halmos

In section \ref{sec:main} we will be mostly interested in applying the above results to constant or periodic inputs $I$ since these are the classical inputs appearing in the studies of Hodgkin-Huxley model.

\section{Application to physiology}\label{sec:main}
In this section we apply the above ideas to a random system based on the Hodgkin-Huxley model well known in physiology. This random system belongs to the family of SDEs with internal variables and random input presented in section \ref{sec:internal}. We start by some reminders on the deterministic Hodgkin-Huxley model that we write (HH) for short. At the end of the following section we also present the general class of conductance-based models of which the Hodgkin-Huxley system is an example as well as the Morris -Lecar model.

\subsection{The deterministic (HH) system}\label{sec:1.0} 

The deterministic Hodgkin-Huxley model for the membrane potential of a neuron (cf \cite{HH-52}) has been extensively studied over the last decades. There seems to be a large agreement that it models adequately many observations made on the response to an external input, in many types of neurons. This model belongs to the family of conductance-based models. Indeed it includes two types of voltage-gated ion channels responsible for the import of Na$^+$ and export of K$^+$ ions through the membrane (for a modern introduction the reader may consult \cite{I-09}).
The Hodgkin-Huxley equations with input $I$ which may be time dependent, is the $4$ dimensional system 
$$
\left\{\begin{array}{l}
dV_t  \;=\; I(t)\, dt \;- \left[\, \ov g_{\rm K}\,n_t^4\, (V_t-E_{\rm K}) \;+\; \ov g_{\rm Na}\,m_t^3\, h_t\, (V_t  -E_{\rm Na}) \;+\; \ov g_{\rm L}\, (V_t-E_{\rm L}) \right] dt\\
dn_t \;=\;  \left[\, \al_n(V_t)\,(1-n_t)  \;-\; \beta_n(V_t)\, n_t  \,\right] dt  \\
dm_t \;=\;  \left[\, \al_m(V_t)\,(1-m_t)  \;-\; \beta_m(V_t)\, m_t  \,\right] dt  \\
dh_t \;=\;  \left[\, \al_h(V_t)\,(1-h_t)  \;-\; \beta_h(V_t)\, h_t  \,\right] dt ,
\end{array}\right. 
\leqno{\rm (HH)}$$
where we adopt the notations and constants of \cite{I-09}.  The parameter $\ov g_{\rm Na}$ (resp. $\ov g_{\rm K}$) is the maximal conductance of a sodium (resp. potassium) channel while $\ov g_{\rm L} $ is the leak conductance. The time dependent conductance of a sodium (resp. potassium) channel depends on the state of four gates which can be open or closed; it is maximal when all gates are open. There are two types of gates $m$ and $h$ for sodium, one type $n$ for potassium. The variables $n_t$, $m_t$, $h_t$ describe the probability that a gate of corresponding type be open at time $t$. The parameter values $\ov g_{\rm K} = 36$, $\ov g_{\rm Na} = 120$, $\ov g_{\rm L} = 0. 3$ , $E_{\rm K} = - 12$, $E_{\rm Na} = 120$, $E_{\rm L} =  10.6$ are those of \cite{I-09}. The parameters $E_{\rm K} $, $E_{\rm Na} $, $E_{\rm L}$ are called reversal potentials. 
The functions $\al_n, \beta_n, \al_m, \beta_m, \al_h, \beta_h$ take values in $(0,\infty)$ and are analytic, i.e. they admit a power series representation on $\bbr$. They are given as follows:
\begin{equation}
\begin{array}{llllll}
\alpha_n(v)  &=& \frac{0.1-0.01v }{\exp(1-0.1v)-1}, &  \beta_n(v) &= &0.125\exp(-v/80) ,  \\
\alpha_m(v)& = &\frac{2.5-0.1v}{\exp(2.5-0.1v)-1} , & \beta_m(v)&= &4\exp(-v/18) ,  \\
\alpha_h(v) &= &0.07\exp(-v/20) , &\beta_h(v) &=& \frac{1}{\exp(3-0.1v)+1}.
\end{array}
\end{equation}
Moreover setting $a_n:=\al_n+\beta_n$, $b_n:=\al_n$ and the analog for $m$ and $h$, the system (HH) fits the notations and assumptions of section \ref{sec:internal}.  If the variable $V$ is kept constant at $v\in\bbr$, the variables $n_t$, $m_t$, $h_t$ converge respectively (cf. Proposition \ref{globallyasymptoticallystable}) towards: 
\begin{equation}\label{eq:ninfty}
n_\infty(v) := \frac{\al_n}{\al_n+\beta_n}(v) \;,\; m_\infty(v) := \frac{\al_m}{\al_m+\beta_m}(v) \;,\; h_\infty(v) := \frac{\al_h}{\al_h+\beta_h}(v) \;.
\end{equation}
In the case of (HH), the function $F$ in section \ref{sec:internal} reads
\begin{equation}\label{eq:F}
F(v,n,m,h) =-[ \ov g_{\rm K}\,n^4\, (v-E_{\rm K})+\ov g_{\rm Na}\,m^3\, h\, (v  -E_{\rm Na}) \;+\ov g_{\rm L}\, (v-E_{\rm L})]. 
\end{equation}

The Hodgkin-Huxley system exhibits a broad range of possible and qualitatively quite different behaviors, depending on the specific input $I$. In response to a periodic input, the solution of (HH) displays a periodic behavior (regular spiking of the neuron on a long time window) only in special situations. Let us first mention that there exists some interval $U$ such that time-constant input in $U$ results in periodic behavior for the solution of (HH) (see \cite{RM-80}). For an oscillating input, there exists some interval $J$ such that oscillating inputs with frequencies in $J$ yield periodic behavior (see \cite{AMI-84}). Periodic behavior includes that the period of the output can be a multiple of the period of the input. However, the input frequency has to be compatible with a range of preferred frequencies of (HH), a fact which is similarly encountered in biological observations (see \cite{I-09}). Indeed there are also intervals $\wt I$ and $\wt J$ such that time-constant input in $\wt I$ or oscillating input at frequency $f\in \wt J$ leads to chaotic behavior. Using numerical methods \cite{Endler} gives a complete tableau.

\subsection{(HH) with random input}

It has been shown in \cite{PTW-10} that conductance-based models like (HH) are fluid limits of a sequence of stochastic models where the states of gates jump from the open state to the closed one with high frequency, when the number of gates of each type goes to infinity. In this framework, the stochastic Hodgkin-Huxley models are Piecewise Deterministic Markov Processes, and limit theorems in this setting enable to study the impact of the {\it channel noise} (also  called {\it intrinsic noise}) on the way the neuron codes the information it receives for instance via spike latency.

Our setting is different. The noise here is external, it comes from the network in which the neuron is embedded, through its dendritic system. This system has a complicated topological structure and carries a large number of synapses which register spike trains emitted from a large number of other neurons within the same active network. The resulting stochastic Hodgkin-Huxley system has the form (\ref{internal}) of a 5-dimensional SDE with internal variables and random input that we recall below:
$$
\left\{\begin{array}{l}
dV_t  \;=\; d \xi_t\, \;- \left[\, \ov g_{\rm K}\,n_t^4\, (V_t-E_{\rm K}) \;+\; \ov g_{\rm Na}\,m_t^3\, h_t\, (V_t  -E_{\rm Na}) \;+\; \ov g_{\rm L}\, (V_t-E_{\rm L}) \right] dt\, ,\\
dn_t \;=\;  \left[\, \al_n(V_t)\,(1-n_t)  \;-\; \beta_n(V_t)\, n_t  \,\right] dt\, ,  \\
dm_t \;=\;  \left[\, \al_m(V_t)\,(1-m_t)  \;-\; \beta_m(V_t)\, m_t  \,\right] dt\, ,  \\
dh_t \;=\;  \left[\, \al_h(V_t)\,(1-h_t)  \;-\; \beta_h(V_t)\, h_t  \,\right] dt \, ,\\
d\xi_t \;=\; b_5(t,\xi_t) dt \;+\sigma (\xi_t)dW_t.
\end{array}\right. 
\leqno{\rm (\xi HH)}
$$
We may also model the cumulated dendritic input as a diffusion of mean-reverting type carrying a deterministic signal $S$. This corresponds to 
$$
\left\{\begin{array}{l}
dV_t  \;=\; d\xi_t \;- \left[\, \ov g_{\rm K}\,n_t^4\, (V_t-E_{\rm K}) \;+\; \ov g_{\rm Na}\,m_t^3\, h_t\, (V_t  -E_{\rm Na}) \;+\; \ov g_{\rm L}\, (V_t-E_{\rm L}) \right] dt\, ,\\
dn_t \;=\;  \left[\, \al_n(V_t)\,(1-n_t)  \;-\; \beta_n(V_t)\, n_t  \,\right] dt \, , \\
dm_t \;=\;  \left[\, \al_m(V_t)\,(1-m_t)  \;-\; \beta_m(V_t)\, m_t  \,\right] dt \, , \\
dh_t \;=\;  \left[\, \al_h(V_t)\,(1-h_t)  \;-\; \beta_h(V_t)\, h_t  \,\right] dt\, ,   \\
d\xi_t \;=\; (\, S(t)-\xi_t\,)\, \tau dt \;+\; \gamma\, q (\xi_t)\, \sqrt{\tau} dW_t \, ,
\end{array}\right.
\leqno{\rm (\xi HHd)}
$$
parametrized in terms of $\tau$ (governing speed) and $\gamma$ (governing spread). For instance $\xi$ can be of Ornstein-Uhlenbeck (OU) type (then $U=\bbr$, $q (\cdot)\equiv 1$) or of Cox-Ingersoll-Ross (CIR) type (then $U=(-K,\infty)$, $q (x)=\sqrt{(x+K)\vee 0\;}$ for $x\in U$, and $K$ is chosen in $]\frac{\gamma^2}{2}+\sup|S|,+\infty[$). Such a choice builds on the statistical study \cite{H-07}. When the deterministic signal $S$ is periodic, it is shown in \cite{HK-10} that $\xi $, either OU type or CIR type, admits a periodically invariant regime under which the signal $S(\cdot) $ is related to expectations of $\xi$ via the formula $s \to E_{\pi, 0} ( \xi_s ) =  \int_0^\infty S(s - \frac{r}{\tau} ) e^{ - r } dr $. Moreover when $\xi$ is of OU type, the periodic ergodicity of the solution to ($\xi$HHd) can be addressed. This is the topic of \cite{HLT-2}. 


\subsection{Weak H\"ormander condition for ($\xi$HH)}\label{secweakhoerHH}
As already noticed, the stochastic Hodgkin-Huxley model $(\xi {\rm HH})$ perfectly fits into Section \ref{ivri}. Therefore, applying Theorem \ref{theo:internal_hoer} and Definition \ref{determinant} we know that we have to consider the $4$-dimensional determinant whose columns are the partial derivatives of the coefficients of (HH) with respect to the first variable, from the first order to the fourth one and look for points where it does not vanish. In the case of (HH), the function $F$, given in (\ref{eq:F}), is linear in the first variable, call it $v$, hence $\partial_v^{(k)}F=0$ for $k\in\{2,3,4\}$. Moreover $\partial_v F(v,n,m,h)=-( \ov g_{\rm K}\,n^4+\ov g_{\rm Na}\,m^3\, h+\ov g_{\rm L})$ which never vanishes on $[0,1]^3$. So actually in this case, it is sufficient to consider a $3$-dimensional determinant extracted from ${\bf D}$. We obtain the

\begin{prop}\label{weakhoerHH}
Assume that $\sigma$ remains positive on $U$. Let us introduce the notation $d_n(v,n) :=  -a_n(v)n+b_n(v)$ and analogous ones for $m$ and $h$. The weak H\"ormander condition {\bf (H4)} for ($\xi$HH) is satisfied at any point $(v,n,m,h,\zeta)\in \bbr\times[0,1]^3\times U$ where ${\Delta}(v,n,m,h)\neq 0$ with 
\begin{equation}\label{eq:det}
{\bf \Delta}(v,n,m,h) \;:=\; \det \left( \begin{array}{lll} 
\partial_v^{(2)} d_n  & \partial_v^{(3)} d_n & \partial_v^{(4)} d_n \\
\partial_v^{(2)} d_m & \partial_v^{(3)} d_m &  \partial_v^{(4)} d_m  \\
\partial_v^{(2)} d_h & \partial_v^{(3)} d_h & \partial_v^{(4)} d_h  \\
\end{array} \right). 
\end{equation}
\end{prop}

\begin{prop}
The set of points in $\bbr\times[0,1]^3\times U$ where the weak H\"ormander condition for ($\xi$HH) is satisfied has full Lebesgue measure.
\end{prop}

{\bf Proof.}
We say that a set has full Lebesgue measure if its complement has Lebesgue measure zero. Firstly it can be shown numerically that indeed there exists points $(v,n,m,h,\zeta)$ where ${\bf \Delta}(v,n,m,h)\neq 0$ (see Section \ref{ex:3} below). Moreover, for any fixed $v \in \bbr ,$ the function $ (n,m,h) \to {\bf \Delta}( v, n,m,h) $ is a polynomial of degree three in the three variables $n,m,h .$ In particular, for any fixed $v,$ either ${\bf \Delta}(v,.,.,.)$ vanishes identically on $(0,1)^3$, or its zeros form a two-dimensional sub-manifold of $(0,1)^3 .$ 
Finally, since ${\bf \Delta}$ is a sum of terms
$$
\mbox{(some power series in $v$)} \cdot n^{\vep_n} m^{\vep_m} h^{\vep_h}
$$
with epsilons taking values 0 or 1, it is impossible to have small open $v$-intervals where it vanishes identically on $(0,1)^3 .$ We conclude the proof by integrating over $v$ and using Fubini's Theorem. \halmos

As pointed out in section \ref{weakHforivri}, Proposition \ref{weakhoerHH} provides a sufficient condition ensuring that the weak H\"ormander condition is satisfied locally. However, this condition is convenient since a numerical study of the involved determinant can be conducted. This is done in section \ref{ex:3} below. Indeed we are not able to characterize the whole set of points where the weak H\"ormander condition holds, unless we make more stringent assumptions on the last component $\xi$, like for instance that the coefficients of its SDE are analytic (cf. \cite{HLT-2}).

\subsection{Numerical study of the determinant ${\bf \Delta}$}\label{ex:3}

We computed numerically the value of ${\bf \Delta}$ at equilibrium points of (HH) associated to different values of constant input $I(t)$ and also along a stable periodic orbit. 

First consider equilibrium points which have the form $(v, n_\infty(v),m_\infty(v), h_\infty(v))$ as in (\ref{eq:ninfty}). The particular point $(0, n_\infty(0), m_\infty(0), h_\infty(0) )$ is the equilibrium point for (HH) with constant input $c$ given by $c := F( 0, n_\infty(0), m_\infty(0), h_\infty(0) ) \approx -0.0534$. We found that ${\bf\Delta}(0, n_\infty(0), m_\infty(0), h_\infty(0) )<0$. Moreover the function $v \mapsto {\bf \Delta}\left(v, n_\infty(v), m_\infty(v), h_\infty(v)\right)$ has exactly two zeros on the interval $I_0=(- 15, + 30)$ located at $v \approx -11.4796$ and $\approx +10.3444$. The function $F_\infty(v):=F(v, n_\infty(v),m_\infty(v), h_\infty(v))$ is strictly increasing on an interval $I$ containing $I_0$. On $I$ there is a bijection between the constant input $I(t)\equiv c$ and the equilibrium value given by the equality $F_\infty(v)$. We used this fact below since it may be more convenient to work with the variable $v$ than with $c$. The corresponding range of values for $c$ is given by $c \in  (F_\infty ( -10) , F_\infty  (+10)) = (-6.15, 26.61) .$ Hence for all values of $c$ belonging to $  (-6.15, 26.61),$ the determinant of the associated equilibrium point remains strictly negative.   

Second we studied $ t \mapsto {\bf \Delta}(v_t,n_t,m_t,h_t)$ along a stable orbit of (HH) with constant input $c=15$. The periodic behavior is shown in Figure 2, starting in a numerical approximation to the equilibrium point which is unstable, and the system switches towards a stable orbit. In this picture, already the last four orbits can be superposed almost perfectly. The value of ${\bf\Delta}$ at equidistant time epochs on the last complete orbit (starting and ending when the membrane potential $v$ up-crosses the level $0$, and having its spike near time $t=180$) is provided in Figure 1. In a window requiring approximately one third of the time needed to run the orbit ${\bf\Delta}(v,n,m,h)$ remains negative and well separated from zero. Very roughly, this segment starts when the variable $v$ up-crosses the level $-2$ and ends when it up-crosses the level $+5$. On the remaining parts of the orbit, ${\bf\Delta}$ changes sign several times. In particular ${\bf \Delta}$ takes values very close to zero immediately after the top of the spike,  i.e.\ after the variable $v$ has reached its maximum over the stable orbit.

\newpage

\begin{figure}[!h]
\begin{center}
   \includegraphics[width=0.80\textwidth]{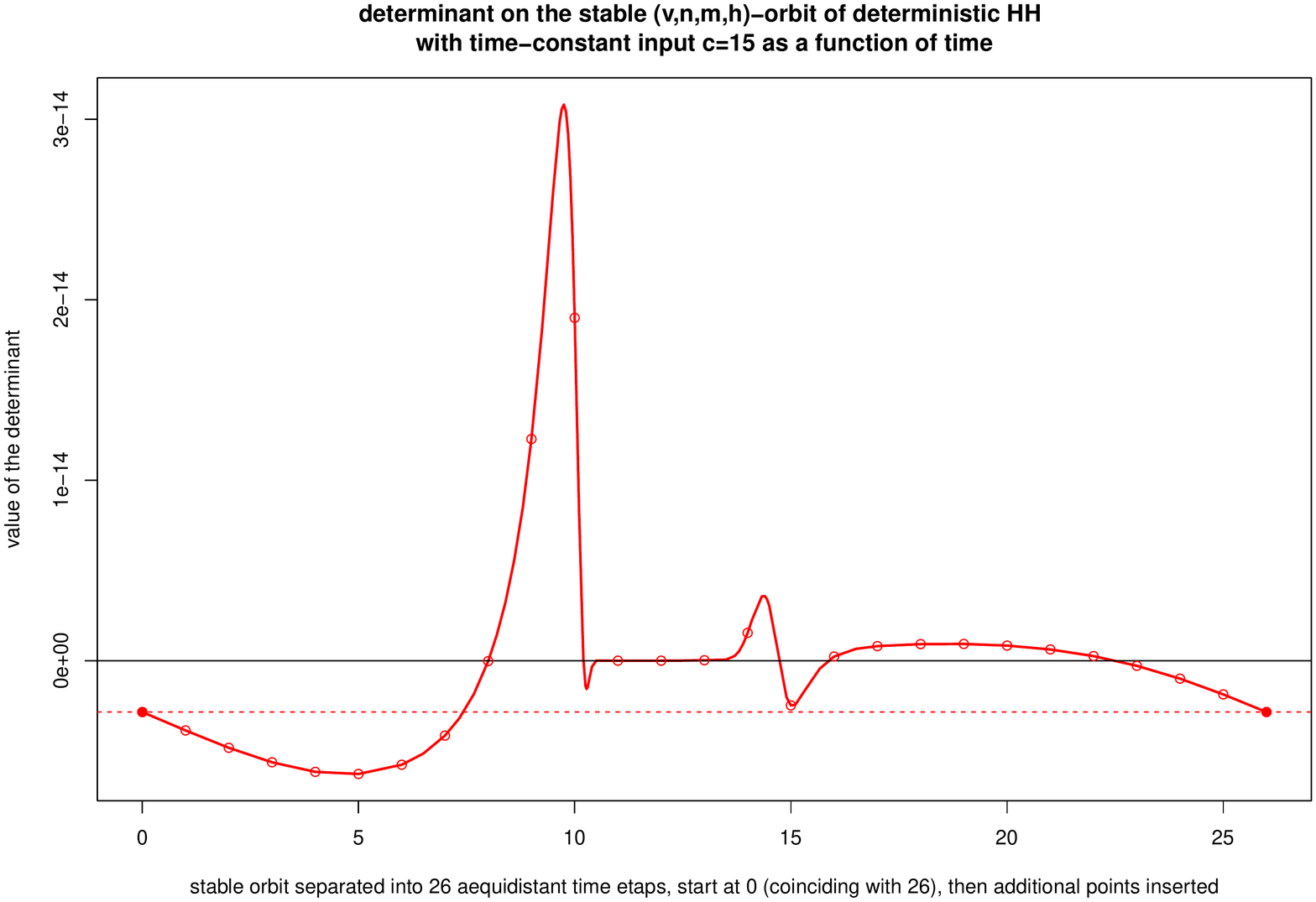} 
\end{center}
\caption{\small Determinant ${\bf \Delta}$ calculated on the stable orbit of the deterministic system (HH) with constant input $c=15$. The time needed to run the orbit is $\approx 12.56$ ms.}
\begin{center}
   \includegraphics[width=0.8\textwidth]{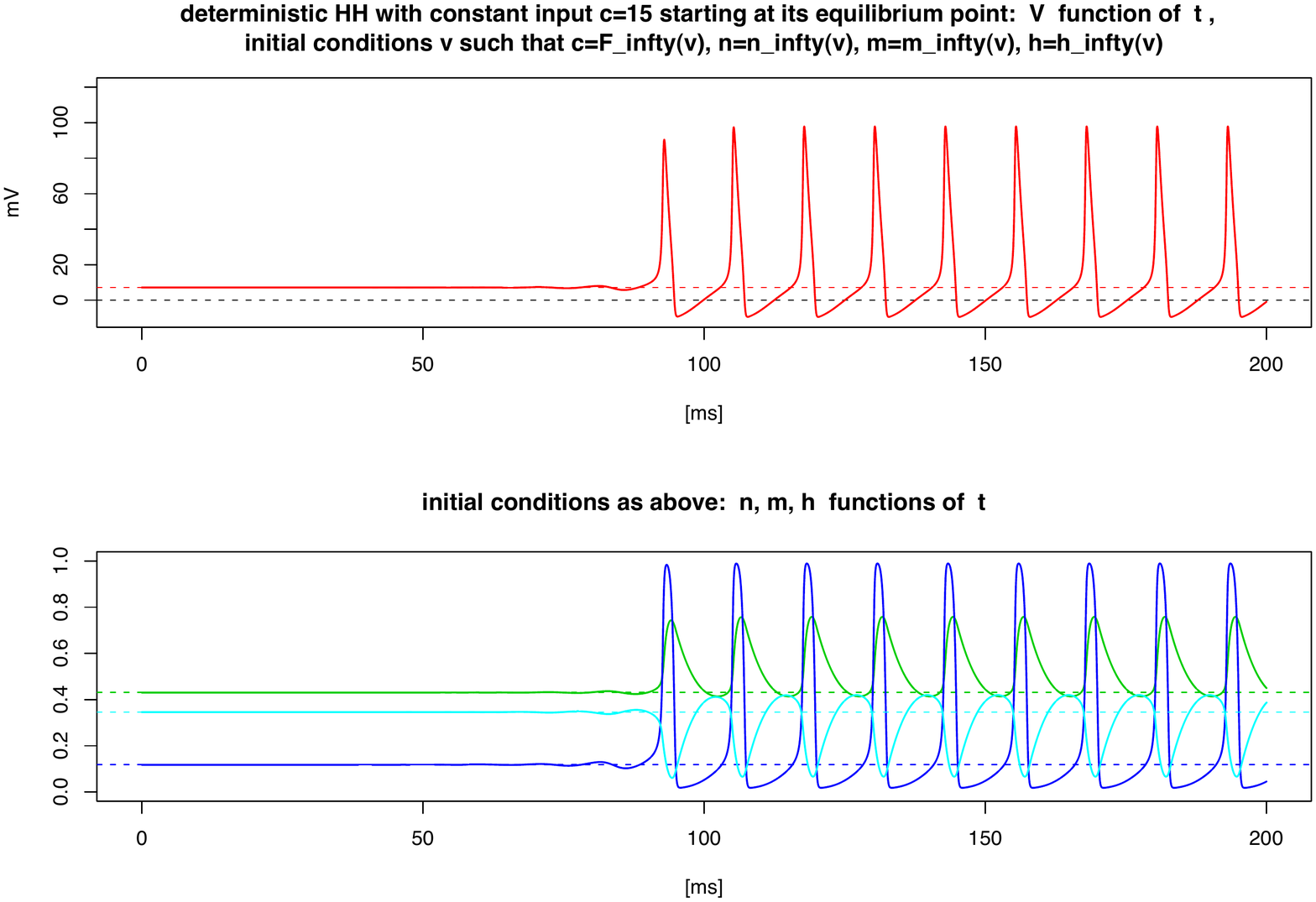} 
\end{center}
\caption{\small Deterministic HH with constant input $c=15.$
}
\end{figure}

\subsection{Positivity regions for ($\xi$HHd)}

We are interested in applying the results of section \ref{positivedensities} mostly in two situations where the inputs are either constant or periodic. For suitable constant input $c$ (see section \ref{ex:3}), fixed $\zeta\in U$ and $t>0$, consider 
\begin{equation}\label{equilibria}
x_c := ( v_c, n_\infty(v_c), m_\infty(v_c), h_\infty(v_c), \zeta ) \quad,\quad x'_c := ( v_c, n_\infty(v_c), m_\infty(v_c), h_\infty(v_c), \zeta+ct )
\end{equation}
where $v_c$ is the unique solution of $F(v_c,n_\infty(v_c), m_\infty(v_c), h_\infty(v_c))=c$. In particular $( v_c, n_\infty(v_c), m_\infty(v_c), h_\infty(v_c) )$ is an equilibrium point of (HH) with constant input $I\equiv c$.

Let us recall that we denote by $P_{s,t}(\cdot,\cdot)_{s<t}$ the semigroup of the process $X = (X_t)_{ t \geq 0} $ which is the solution of the stochastic system ($\xi$HHd) driven by the fixed signal $S.$  

\begin{prop}\label{cor:1} Assume that $\zeta+cs\in U$ for all $0\le s\le t$. Consider $x_c$ and $x'_c$ defined in (\ref{equilibria}). Then for all $\vep>0$, there exists $\delta > 0$ such that for all $ x'' \in B_\delta ( x_c) , $ $    P_{ 0, t} ( x'' , B_{\varepsilon }( x_c') ) > 0$.
\end{prop} 

\begin{prop}\label{cor:2} 
We keep the assumptions and notations of Proposition \ref{cor:1} and we assume moreover that ${\bf \Delta}( v_c, n_\infty(v_c), m_\infty(v_c), h_\infty(v_c) )\neq 0$. There exists $\delta>0$ such that 
for $ K_c = \overline B_\delta ( x_c)   $ and $ K_c' = \overline B_\delta (x_c') ,$  
$$ \inf_{ x \in K_c } \inf_{x' \in K_c'}  p_{ 0, t} (x, x') > 0 $$.
\end{prop}
Remember that the assumption ${\bf \Delta}( v_c, n_\infty(v_c), m_\infty(v_c), h_\infty(v_c) )\neq 0$ ensures that the local H\"ormander condition holds at both points $x_c$ and $x'_c$. We have checked numerically that this assumption is satisfied for $c\in ]-6.15, 26.61[$ (cf. section \ref{ex:3}).

The second situation which we consider is the deterministic system (HH) receiving a sinusoidal signal 
$
I(t)  = a \left( 1 + sin( 2\pi \frac{t}{T} ) \right) , 
$
parametrized by $(a,T)$ where $a>0$ is some constant. This system presents additional features (see \cite{AMI-84} for a modified system, see \cite{Endler} for the above (HH) system). There are specified subsets $D_1$, $D_2$, $D_3$, $D_4$ in $(0,\infty)\times(0,\infty)$ with the following properties: for $(a,T)$ in $D_1$ (HH) is periodic with small oscillations which cannot be interpreted as spiking. For  $(a,T)$ in $D_2$ the system moves on a $T$-periodic orbit, and the projection $t\mapsto V_t$ resembles the membrane potential of a regularly spiking neuron (single spikes or spike bursts per orbit). For  $(a,T)$ in $D_3$ the system is periodic with period a multiple of $T$. For $(a,T)$ in $D_4$ it behaves irregularly and does not exhibit periodic behavior. 

For this particular input $I$ with parameters $(a,T)\in D_2$, consider the points
$$
x:= (0, n^*, m^*, h^*, \zeta ) \quad,\quad x':= (0, n^*, m^*, h^*, \zeta + \int_0^T I(r) dr )
$$ 
such that $(0, n^*, m^*, h^*)$ corresponds to exactly one point on the stable orbit of (HH) at which the membrane potential equals $0$. Since we choose $(a,T)\in D_2$, the solution of (HH) performs exactly one tour on the stable orbit during $[0,T]$. 
\begin{prop}\label{cor:3}
Assume that $\zeta+\int_0^t I(r) dr \in U$ for all $0\le t\le T$ and that ${\bf \Delta}(x)\neq 0$. There exists
$ \delta > 0 $ such that for $ K  = \overline B_\delta ( x)   $ and $ K' = \overline B_\delta (x') ,$  
$$ \inf_{ y \in K  } \inf_{y' \in K' }  p_{ 0, T} (y, y') > 0 .$$
\end{prop}
Note that assuming ${\bf \Delta}( 0, n^*, m^*, h^* )\neq 0$ implies that the local H\"ormander condition holds at both points $x^*$ and $z^*$. In our example of Figure 1 this is satisfied at the point on the orbit at which the membrane potential \ul{up-crosses} level $0$. 

\subsection{Conductance-based models.}
As previously mentioned the Hodgkin-Huxley model is a conductance-based model. Such models are built to reflect the complex electro-chemical process at work at the level of the cell membrane. In these models the variable of interest is the difference of potential between the inside and the outside of the cell often called {\it membrane potential}. The time evolution of this potential is strongly coupled with the state of units which are channels located in the membrane. They can also be components of these channels called gates like in the Hodgkin-Huxley system. When they are open they allow selectively specific ions to flow inside or outside the membrane creating gradients of charges.  The equation of the potential is given by a current balance (Kirchoff's Law). In these models the current created by the ionic flow in a channel permeable to species $i$ ions obeys Ohm's Law and has the form $G_i (V-{\overline V}_i)$ where ${\overline V}_i$ is constant. The coefficient $G_i$ is called the {\it conductance}. Finally, the repartition of channels is not homogeneous across the membrane and there are parts of it which do not contain channels. These parts contribute by a current ${\overline G}_L(V-{\overline V}_L)$ where the conductance ${\overline G}_L$ is constant (so is ${\overline V}_L$). Then the equation for $V_t$ has the form
\begin{equation}\label{conductance-based1}
{\overline C}dV_t=-[\sum_{i=1}^k G_i(V_t-{\overline V}_i)+{\overline G}_L(V-{\overline V}_L)]\, ,
\end{equation} 
where ${\overline C}$ is the membrane capacitance which is constant and $k$ denotes the number of ion species involved.\\ 
The conductance $G_i$ depends on the probability that a channel permeable to species $i$ ions is open. Let us denote by $p_{i,t}$ the probability that such a channel is open at time $t$.  Then the general form of a conductance-based model reads as follows
\begin{eqnarray}\label{conductance-based2}
{\overline C}dV_t&=&-[\sum_{i=1}^k G_i(p_i\, , i \leq 1\leq k)(V_t-{\overline V}_i)+{\overline G}_L(V-V_L)]dt\, ,  \\
dp_{i,t}&=&[-a_i(V_t)p_{i,t}+b_i(V_t)]dt, \, \, \, i = 1, \cdot, \cdot, \cdot , k.\nonumber
\end{eqnarray}
These systems are fully coupled which makes their study complex. However it is important to note that the right hand-side of the $V$-equation as a function of the vector $(v, p_i\, , i \leq 1\leq k)$ is linear w.r.t. $v$. We will use this fact in Proposition \ref{weakhoerCB} below. It is easy to see that these systems belong to the class of systems with internal variables of the form (\ref{subinternal}). In particular they can be represented as fluid limits of sequences of Piecewise Deterministic Markov Processes as explained in section \ref{sec:internal}.\\
In the Hodgkin-Huxley system ${\rm (HH)}$ the units of interest are gates not channels. A potassium channel has four gates of the same type (the classical notation is $n$). A sodium channel is made of four gates of two types (classical notation is $m$ and $h$). The conductances are then functions of the potential and the probabilities that a gate of a given type be open at time $t$. From their data, Hodgkin and Huxley (see \cite{HH-52}) proposed $G_ {\rm Na}=\overline{g}_{\rm Na}m_t^3 h_t$ for a sodium channel and $G_ {\rm K}=\overline{g}_{\rm K}n_t^4$ for a potassium channel, which turned out to be a far reaching choice since at that time it was not possible experimentally to have access to the channels. The gates structure was confirmed many years later.

Stochastic conductance-based models can be also considered when systems like (\ref{conductance-based2}) are submitted to a random input. From (\ref{conductance-based2}) on obtains the $k+2$-dimensional system
\begin{eqnarray}\label{conductance-based3}
{\overline C}dV_t&=&d\xi_t -[\sum_{i=1}^k G_i(p_i\, , i \leq 1\leq k)(V_t-{\overline V}_i)+{\overline G}_L(V-V_L)]dt\, ,  \\
dp_{i,t}&=&[-a_i(V_t)p_{i,t}+b_i(V_t)]dt, \, \, \, i = 1, \cdot, \cdot, \cdot , k\, ,\nonumber\\
d\xi_t &=& b_{k+2}(t,\xi_t) dt \;+\sigma (\xi_t)dW_t.\nonumber
\end{eqnarray}
for a general random input $(\xi_t)$ or the following one when $(\xi_t)$ carries a deterministic signal
\begin{eqnarray}\label{conductance-based4}
{\overline C}dV_t&=&d\xi_t -[\sum_{i=1}^k G_i(p_i\, , i \leq 1\leq k)(V_t-{\overline V}_i)+{\overline G}_L(V-V_L)]dt\, ,  \\
dp_{i,t}&=&[-a_i(V_t)p_{i,t}+b_i(V_t)]dt, \, \, \, i = 1, \cdot, \cdot, \cdot , k\, ,\nonumber\\
d\xi_t &=& (\, S(t)-\xi_t\,)\, \tau dt \;+\; \gamma\, q (\xi_t)\, \sqrt{\tau} dW_t.\nonumber
\end{eqnarray}

We now state for system (\ref{conductance-based3}) the following sufficient condition which implies that the weak H\"ormander condition is satisfied. 
\begin{prop}\label{weakhoerCB}
Assume that $\sigma$ remains positive on $U$. For $1\leq i\leq k\, $, let us set $d_i(v,p_i) :=  -a_i(v)p_i+b_i(v)$. Let ${\Delta}(v, p_i\, , 1\leq i\leq k)$ be the $k$-dimensional determinant whose lines are the partial derivatives of the $d_i$ w.r.t. $v$, from order $2$ to order $k+1$. Consider a point $(v,p_i\, , 1\leq i\leq k\, , \zeta)$ such that $\sum_{i=1}^k G_i(p_i\, , i \leq 1\leq k)+{\overline G}_L\neq 0$.  Then the weak H\"ormander condition {\bf (H4)} for (\ref{conductance-based3}) holds at this point if it satisfies ${\Delta}(v, p_i\, , 1\leq i\leq k)\neq 0$.
\end{prop}

\noindent {\bf Proof of Proposition \ref{weakhoerCB}.} The proof of this proposition is similar to Proposition \ref{weakhoerHH}. We use the linearity of the r.h.s. of (\ref{conductance-based1}) w.r.t. the variable $v$ and the assumption $\sum_{i=1}^k G_i(p_i\, , i \leq 1\leq k)+{\overline G}_L\neq 0$. This latter condition is always satisfied in the case of ${\rm (\xi HH)}$ since the conductances are positive and do not vanish simultaneously.\halmos

\section{Appendix : Proof of (\ref{eq:classical})}
In this appendix section we prove (\ref{eq:classical}). To do this, we have to extend some basic arguments from Malliavin Calculus to the non time homogeneous case. For our strongly degenerate diffusion, we have to work with a local H\"ormander condition of order $L$, using a localization argument. In \cite{stefano} local ellipticity is supposed to hold, hence the order is $L=1$.  Thus our local H\"ormander condition of order $L$ changes the order of growth in some bounds ( see also \cite{Fournier} and \cite{Bally}). For the basic concepts of Malliavin calculus, we refer the reader to the classical reference \cite{Nualart}.

Throughout this section we fix $T>0$ and denote by $\bar X  $ the unique strong solution of the SDE (recall also (\ref{eq:processgood}))
\begin{equation}\label{eq:processgoodbis}
\bar X_{i,t} = x_i + \int_0^t \bar b_i ( s, \bar X_s) ds + \int_0^t \bar \sigma_i (\bar X_s) d W_s , \, t \le T ,\,  1 \le i \le m ,
\end{equation}
where $x \in \bbr^m ,$ $ \bar \sigma_i \in C^\infty_b ( \bbr^m ) ,  $  and where $\bar b_i (t,x) $ and all partial derivatives $ \partial^\alpha_x \partial^\beta_t \bar b_i (t, \cdot ) $ are bounded uniformly in $t \in [0, T ]$.
Let 
$$ \tilde{\bar b}_i (t,x) = \bar b_i (t,x) - \frac12 \sum_{k=1}^m \bar\sigma_k (x) \frac{ \partial \bar \sigma_i }{\partial x_k } (x) ,\;  1 \le i \le m ,$$
and $\bar A_0  = \frac{\partial}{\partial t } + \tilde{\bar b}  $ be the associated time-space directional derivative. We use notation analogous to section \ref{sec:hoerm} and put $ \bar A_1 = \bar \sigma .$
The local H\"ormander condition at a point $x$ for a given number of brackets $L$ is expressed by $  \bar {\cal V}_L (x) > 0$ where we set $ \bar {\cal V}_L (x):= \inf_{ 0 \le t \le T }  {\cal V}_L (t, x) $ (see section \ref{sec:hoerm}).  

The main ingredient for the control of the weight in Malliavin's integration by parts formula as in formula (\ref{eq:classical}) is to obtain estimates of Malliavin's covariance matrix. We check that all results obtained in \cite{KS} are still valid in our framework. Let 
$$ (Y_t)_{i,j} = \frac{ \partial \bar X_{i,t}}{\partial x_j} , 1 \le i, j \le m .$$
Then $Y$ satisfies the following linear equation having bounded coefficients (bounded with respect to time and space)
$$ Y_t = I_m + \int_0^t  \partial \bar b ( s, \bar X_s) Y_s  ds + \int _0^t \partial \bar \sigma ( \bar X_s) Y_s d W_s .$$
Here $I_m$ is the $m\times m-$unity matrix and $ \partial \bar b $ and $ \partial \bar \sigma $ are the $m\times m-$matrices having components $ ( \partial \bar b)_{i,j} (t, x) = \frac{\partial \bar b_i }{\partial x_j} (t, x ) $ and  $ ( \partial \bar \sigma)_{i,j} ( x) = \frac{\partial \bar \sigma_i }{\partial x_j} ( x )  .$
By means of It\^o's formula, one shows that $Y_t$ is invertible. The inverse $Z_t$ satisfies the linear equation with coefficients bounded in $t$ and in $x$ given by 
\begin{equation}\label{eq:z}
Z_t = I_m - \int_0^t \partial \tilde{\bar b} (s,\bar X_s)  Z_s ds - \int_0^t  \partial \bar \sigma ( \bar X_s) Z_s \circ dW_s   ,
\end{equation}
where $ \circ d W_s $ denotes the Stratonovitch integral. In this framework, the following estimates are classical (see e.g.  \cite{KS} or \cite{stefano}). 
For all $ 0 \le s \le t \le T ,$ for all $ p \geq 1,$ 
\begin{equation}\label{eq:ub1}
E\left( \sup_{ r :s \le r \le t } | \bar X_{i,r} - \bar X_{i,s} |^p \right) \le C(T,p,m, \bar b , \bar \sigma ) (t-s)^{p/2} ,
\end{equation}
\begin{equation}\label{eq:ub2}
\sup_{ s \le t } E ( |(Z_s)_{i,j}|^p) \le C(T,p,m, \bar b , \bar \sigma ) ,\; 1 \le i, j \le m ,
\end{equation} 
\begin{equation}\label{eq:ub3}
\sup_{ r_1 , \ldots , r_k \le t } E \left(  | D_{r_1, \ldots , r_k} \bar X_{i,t}|^p  \right) \le C(T,p,m,k, \bar b , \bar \sigma )  \left( t^{1/2}  +  1\right)^{ (k+1)^2 p} ,
\end{equation}
where the constants depend only on the bounds of the space derivatives of $\bar b$ and $\bar \sigma .$ Notice that the above estimates are not sharp, and much better estimates can be obtained. However, for our purpose, the above estimates are completely sufficient.

As indicated before, the main issue in order to prove (\ref{eq:classical}) is to obtain estimates on the Malliavin covariance matrix. 
So let $ (\sigma_{ \bar X_t})_{i,j}  = < D \bar X_{i,t} , D \bar X_{j,t}>_{L^2 [0, t ]} , 1 \le i, j \le m .$
Then it is well known (see for example \cite{Nualart}, page 110, formula (240)) that 
$$ \sigma_{ \bar X_t} = Y_t \left( \int_0^t Z_s \bar \sigma ( \bar X_s) \bar \sigma^* ( \bar X_s) Z_s^* ds \right) Y_t^* .$$
In order to evaluate the integral, one has to control expressions of the form 
$ Z_s V(s, \XXX_s) ,$ where $ V (t,x) $ is a smooth function of $t$ and $x.$ Using partial integration we obtain (cf. \cite{KS}, formula (2.10))
\begin{multline}\label{eq:important2}
Z_t V (t, \XXX_t ) = V(0, x) + \int_0^t Z_s [ \bar \sigma , V ] ( s, \XXX_s) \circ d W_s 
+ \int_0^t Z_s  [ \frac{\partial}{\partial t} + \tilde{ \bar b} , V] (s, \XXX_s) ds \\
= V(0, x) + \int_0^t Z_s [ \bar \sigma , V ] ( s, \XXX_s) \circ d W_s 
+ \int_0^t Z_s  [ \bar V_0 , V] (s, \XXX_s) ds 
\end{multline}

Iterating (\ref{eq:important2}) we obtain for any $ L \geq 1 , $
\begin{equation}\label{eq:important1}
Z_s \bar \sigma (\XXX_s) = \sum_{ \alpha : \| \alpha \| \le L- 1 } W^{(\alpha)} (s) (\bar A_1)_{(\alpha)} (0, x ) + R_L (s, x, \bar A_1 ),
\end{equation}
where $R_L$ is a remainder term and where $ W^{(\alpha)} $ is a multiple Wiener integral. This is completely analogous to Theorem 2.12 of \cite{KS}. Here, $\bar A_1 = \bar \sigma $ and the $ (\bar A_1)_{(\alpha)} (0, x )$ are the successive Lie brackets. The most important feature in the above development (\ref{eq:important1}) is that the behavior of the remainder term depends only on the supremum norms of derivatives with respect to time and space of $\bar b $ and with respect to space of $\bar \sigma .$ Then, following \cite{KS}, we obtain the

\begin{cor}[Corollary 3.25 of  \cite{KS}]
For any  $p \geq 1 $ and $ t \le 1 , $ for any $ L \geq 1, $ for any $x$ such that $\bar {\cal V}_L (x)>0\, $,
\begin{equation}\label{eq:important3}
E_x \left( | det \sigma_{\XXX_t}|^{-p} \right)^{1/p} \le C(p,m, L) \frac{1}{(\bar {\cal V}_L (x)^{ 1 + \frac{2}{L}} t )^{mL}} .
\end{equation}
\end{cor}

This bound (\ref{eq:important3}) plays the role of the bound (2.20) in \cite{stefano} where for $t$ close to zero, the bound (2.20) is of order $t^{-m}$ due to the local ellipticity condition, while our bound is of order  $t^{-mL}$ due to our condition $\bar {\cal V}_L (x)>0$. Our (\ref{eq:ub3}) is the same bound as (2.17) in \cite{stefano}. Inserting our upper bounds (\ref{eq:important3}) and (\ref{eq:ub3}) in formula (2.25) in the proof of Theorem 2.3 of in \cite{stefano} replaces the r.h.s. obtained there by $O\left( t^{-p\left\{ m[L-1] + 1 \right\}} \right)$ for small $t$. 
With such changes, the argument developed in \cite{stefano} goes through, and we end up in our case with a r.h.s.\ $ O\left( t^{-m\, n_L } \right) $ for small $t$ with some positive constant $n_L$ depending on the order $L$ in $\bar {\cal V}_L (x)$. This allows to obtain the desired formula (\ref{eq:classical}): the factor $\delta^{ - m n_L}$ comes from the fact that we apply the Malliavin calculus over a time interval of total length $\delta$.  Here, since we do not need them, we do not take into account the precise form of the constants reported in \cite{stefano}. 

\section*{Acknowledgments}
We thank Vlad Bally and Michel Bena\"im for very stimulating discussions.

\end{document}